\documentclass[11pt, onesided, leqno]{amsart}
\usepackage{enumerate}
\usepackage[dvips]{graphics}
\usepackage{amsmath,amssymb,mathrsfs}
\newcommand{\doublespace}
   {\addtolength{\baselineskip}{0.15\baselineskip}}

\newtheorem{pdef}{Definition}[section]
\newtheorem{thm}[pdef]{Theorem}
\newtheorem{cor}[pdef]{Corollary}
\newtheorem{cond}[pdef]{Condition}
\newtheorem{lem}[pdef]{Lemma}
\newtheorem{exam}[pdef]{Example}
\newtheorem{remark}[pdef]{Remark}
\newtheorem{prop}[pdef]{Proposition}

\newcounter{equationnumber}
\renewcommand{\theequation}{\thesection.\arabic{equation}}
\def\maxconv{\square\kern-.72em\lower-.14ex\hbox{$\vee$}}
\def\maxconvexp{\square\kern-.572em\lower-.025ex\hbox{\footnotesize$\vee$}}
\def\maxconvlow{\square\kern-.576em\lower-.03ex\hbox{\footnotesize$\vee$}}
\def\mathletters{
    \addtocounter{equation}{1}
    \edef\@currentlabel{\theequation}
    \setcounter{equationnumber}{\value{equation}}
    \setcounter{equation}{0}
    \edef\theequation{\@currentlabel\noexpand\alph{equation}}
    }

\title{Bi-Free Extreme Values}

\author{Hao-Wei Huang and Jiun-Chau Wang}

\address{Department of Applied Mathematics, National Sun Yat-sen University,
No. 70, Lienhai Road, Kaohsiung 80424, Taiwan, R.O.C.}

\email{hwhuang@math.nsysu.edu.tw}

\address{Department of Mathematics and Statistics, University of
Saskatchewan, Saskatoon, Saskatchewan S7N 5E6, Canada}
\email{jcwang@math.usask.ca}

\date {\today}

\begin{document}
\maketitle \doublespace \pagestyle{myheadings} \thispagestyle{plain}
\markboth{   }{ }

\begin{abstract} In this paper, we continue Voiculescu's recent work on the analogous extreme value theory in the context of bi-free probability theory. We derive various equivalent conditions for a bivariate distribution function to be bi-freely max-infinitely divisible. A bi-freely max-infinitely divisible distribution function can be expressed in terms of its marginals and a special form of copulas. Such a distribution function is shown to be also max-infinitely divisible in the classical sense. In addition, we characterize the set of bi-free extreme value distribution functions. A distribution function of this type is also bi-free max-stable and represented by its marginals and one copula composing of a Pickands dependence function, as in the classical extreme value theory. As a consequence, the determination of its bi-free domain of attraction is the same as the criteria in the classical theory. To illustrate these connections, some concrete examples are provided.
\end{abstract}


\section{Introduction} The main aim of this paper is to further investigate the introductory of the bi-free extreme value theory lately laid down by Voiculescu \cite{V15EXT}.

In classical probability theory, extreme value distributions arise as limiting distributions for the componentwise maximums or minimums (extreme values) of a sample of independent and identically distributed multivariate random vectors, as the sample size increases. Since the introduction by Fisher and Tippett \cite{EVT}, there have been extensive studies in this subject whether in the univariate or multivariate situations in the past literature. There are only three types of possible limit laws in the univariate case, which are respectively known as the Gumbel, Frechet, and Weibull distributions. A wide variety of situations involving extreme events are inherently multivariate. Various alternative ways have been proposed to describe the dependence structure of a multivariate extreme value distribution function. A bivariate extreme value distribution function, for instance, is determined by its marginals and its unique Pickands dependence function.

Apart from the studies of max-domain of attraction and max-stability, the analysis of the max-infinite divisibility of distribution functions also plays an essential role in dealing with multivariate extremal processes. Every univariate distribution function is max-infinitely divisible, whereas this is no longer true in the multivariate case. For example, the necessary and sufficient conditions for the bivariate normal distribution to be max-infinitely divisible is that its correlation coefficient is non-negative. In general, a multivariate distribution function is max-infinitely divisible if and only if the negative of its logarithm is an exponent measure \cite{Maxclassical}.

It was recently demonstrated in \cite{V15EXT} that for two bi-free bipartite two-faced pairs of self-adjoint operators in a $C^*$-probability space, one can compute the joint distribution of their componentwise maximums. Here the maximum operation on self-adjoint operators is taken with respect to the spectral order relation \cite{Ando}, while the notion of bi-freeness in bi-free probability theory substitutes for the independence in classical probability theory \cite{V14}. This max-operation in free probability theory gives rise to the free max-convolution of univariate distribution functions, and there is a simple formula for computing this convolution \cite{freeextre}. As in extreme value theory, there are three possible limit distributions in free extreme value analysis, governed by the generalized Pareto laws, and each of them has its correspondent in the classical situation. Classical and free max-domains of attraction of the corresponding laws turn out to be the same. In the bi-free framework, the componentwise max-operation induces the bi-free max-convolution of distribution functions on the plane, and a concise formula for dealing with this convolution also exists. Passing to bivariate distribution functions, this binary operation transforms all questions regarding bi-free extremes into problems in classical analysis.

Following the pioneering work of Voiculescu, we would like to contribute to the research of characterizing the set of bivariate distribution functions having the properties of bi-free max-stability and bi-freely max-infinite divisibility. In probability theory and statistics, copulas have become a powerful tool as they provide an elegant and detailed description of the dependence structure among random variables. In the present paper, a particular form of copulas is discovered to describe how the marginals of a bi-freely max-infinitely divisible distribution function are linked. Copulas which are used for the description and of great interest are those in the Ali-Mikhail-Haq family, also one of the most prominent bivariate Archimedean copulas. Thanks to a newly defined transform of bivariate distribution functions in the paper, we show that a bi-freely max-infinitely divisible distribution function under this transform serves as an exponent measure in the theory of max-infinite divisibility. As a matter of fact, any exponent measure having support bounded from below arises in the way described above.

We offer thorough descriptions of bi-free extreme value distribution functions, including their bi-free max-domains of attraction, and the characterization of bi-free max-stability. Specifically speaking, we show that a distribution function is a bi-free extreme value distribution function if and only if it can be represented by means of its marginals and a copula involving a unique Pickands dependence function. Also, a bi-free extreme value distribution function is bi-free max-stable and vice versa. We further derive that a bivariate extreme value distribution function and a bi-free extreme value one, whose marginals are the corresponding classical and free extreme types, determined by a common Pickands dependence function have an identical max-domain of attraction. These findings establish rigid conjunctions between bi-free and classical extreme value theory. To illustrate these established theories, we present a couple of examples including parametric families of copulas widely used in practice.

The organization of the paper is as follows. After the basics of extreme value theory and copulas in Section 2, we derive the limit theory for bi-free max-convolution in Section 3. Section 4 is devoted to studying the bi-freely max-infinite divisibility. A representation and examples of bi-freely max-infinitely divisible distribution functions are provided in Section 5 and 6. The connection between bi-freely and classical max-infinite divisibility is studied in the following section. Section 8 focuses on bi-free extreme value distribution functions and their max-domains of attraction.

\section{Preliminary} Throughout the paper, the bold letter $\mathbf{x}$ is used to
denote a constant vector $(x_1,x_2)$ in $\mathbb{R}^2$. We will also denote by $\boldsymbol\infty=(\infty,\infty)$ and $-\boldsymbol\infty=(-\infty,-\infty)$ in the extended plane. If not otherwise stated, all expressions such as $\mathbf{x}\vee\mathbf{y}=\max\{\mathbf{x},\mathbf{y}\}$ and $\mathbf{x}\leq\mathbf{y}$ are intended to be componentwise operations; for instance, $\mathbf{x}\leq\mathbf{y}$ means $x_j\leq y_j$ for $j=1,2$. Given two points $\mathbf{x}\leq\mathbf{y}$ in $\mathbb{R}^2$, a bounded rectangle $[\mathbf{x},\mathbf{y}]=\{\mathbf{z}:x_j\leq z_j\leq y_j\;\;\mathrm{for}\;j=1,2\}$ is the generalized closed interval in $\mathbb{R}^2$. Analogous notations also apply to rectangles $(\mathbf{x},\mathbf{y})$ and $[\mathbf{x},\mathbf{y})$.

For a real-valued function $H$ defined on some convex set $\mathcal{S}$ in $\mathbb{R}^2$, the $H$-volume of a bounded rectangle $[\mathbf{x},\mathbf{y}]$ contained in $\mathcal{S}$ is defined by
\[V_H\big([\mathbf{x},\mathbf{y}]\big)=H(\mathbf{y})-H((x_1,y_2))-H((y_1,x_2))+H(\mathbf{x}).\] Then $H$ is said to be \emph{quasi-monotone} on $\mathcal{S}$ if $V_H([\mathbf{x},\mathbf{y}])\geq0$ for any $[\mathbf{x},\mathbf{y}]\subset\mathcal{S}$

\begin{remark} \label{qmonotone}
\emph{The quasi-monotonicity of $H$ on a convex set $\mathcal{S}$ ensures that both functions $H(y_1,\cdot)-H(x_1,\cdot)$ and $H(\cdot,y_2)-H(\cdot,x_2)$ increase on their domain whenever points $x_1<y_1$ and $x_2<y_2$ are fixed. Particularly, for $\mathcal{S}=\mathbb{R}^2$ or $\mathcal{S}$ of the form $[\mathbf{L},\boldsymbol\infty)$ or $(\mathbf{L},\boldsymbol\infty)$ with $\mathbf{L}\in\mathbb{R}^2$, we come to the conclusion
\begin{equation} \label{quasiineq}
H(\mathbf{y})-H(\mathbf{x})\leq\lim_{\alpha,\beta\to\infty}\big[H(\alpha,y_2)-H(\alpha,x_2)+H(y_1,\beta)-H(x_1,\beta)\big]
\end{equation} for any points $\mathbf{x}\leq\mathbf{y}$ in $\mathcal{S}$. Another convex set that is generally considered herein is the unit square $[0,1]^2$ or $(0,1]^2$, in which case the same inequality also holds except the points $\alpha$ and $\beta$ in the limit (\ref{quasiineq}) increase to one.}
\end{remark}

\subsection{Copula} In the following, we outline some basics in the bivariate analysis we shall need in what follows. For a random vector $\mathbf{X}=(X_1,X_2)$ on $\mathbb{R}^2$, the distribution function
\[F(\mathbf{x})=\mathrm{Pr}[\mathbf{X}\leq\mathbf{x}],\;\;\;\;\;\mathbf{x}\in\mathbb{R}^2,\] provides a detailed description of dependence structure between the random variables $X_1$ and $X_2$. Another alternative for this is the distribution function, defined as
\[\mu(B)=\mathrm{Pr}[\mathbf{X}\in B]\] for any Borel set $B$ in $\mathbb{R}^2$.
The distribution functions of $X_1$ and $X_2$ are referred to as marginal distribution functions of $F$, denoted by $F_1$ and $F_2$, respectively. Marginal distributions $\mu_1$ and $\mu_2$ of $\mu$ are defined in the same manner.

A distribution function $F$ in two variables apparently enjoys the properties (\emph{i}) $\lim_{x_j\to-\infty}F(\mathbf{x})=0$ for $j=1,2$; (\emph{ii}) $\lim_{\mathbf{x}\to\boldsymbol\infty}F(\mathbf{x})=1$; (\emph{iii}) $F$ is \emph{continuous from above} in the sense that $F(\mathbf{y})\to F(\mathbf{x})$ whenever $\mathbf{y}\to\mathbf{x}$ with $\mathbf{y}\geq\mathbf{x}$; (\emph{iv}) $F$ is \emph{increasing}, i.e., $F(\mathbf{x})\leq F(\mathbf{y})$ for any points $\mathbf{x}\leq\mathbf{y}$ in $\mathbb{R}^2$; and (\emph{v}) $F$ is quasi-monotone on $\mathbb{R}^2$.

Conversely, any function $F:\mathbb{R}^2\to[0,1]$ satisfying the previously mentioned properties (\emph{i})-(\emph{v}) is associated with a unique Borel probability measure $\mu$ on $\mathbb{R}^2$ so that
\[F(\mathbf{x})=\mu((-\boldsymbol\infty,\mathbf{x}])\] holds for any point $\mathbf{x}\in\mathbb{R}^2$ \cite{Billingsley}.

Recall that the \emph{tail distributions} of $F_j$ and $F$ are defined as
\[\overline{F_j}(x)=1-F_j(x)\;\;\;\;\;\mathrm{and}\;\;\;\;\;
\overline{F}(\mathbf{x})=1+F(\mathbf{x})-F_1(x_1)-F_2(x_2)\] for any $x\in\mathbb{R}$ and $\mathbf{x}\in\mathbb{R}^2$. Obviously, we have $\overline{F_j}(x)=\mu_j((x,\infty))$ and $\overline{F}(\mathbf{x})=\mu((\mathbf{x},\boldsymbol\infty))$ if $\mu$ is the Borel probability measure having $F$ as the distribution function.

Next, applying (\ref{quasiineq}) to a bivariate distribution function $F$ results in the estimate
\begin{equation} \label{upperbound}
F(\mathbf{y})-F(\mathbf{x})\leq F_1(y_1)-F_1(x_1)+F_2(y_2)-F_2(x_2),\;\;\;\;\;
\mathbf{x}\leq\mathbf{y}\in\mathbb{R}^2.
\end{equation} An immediate consequence of (\ref{upperbound}) is that the set $\mathcal{C}(F)$ of continuities of $F$ contains $\mathcal{C}(F_1)\times\mathcal{C}(F_2)$, the Cartesian product of the sets of continuities of $F_1$ and $F_2$.

We now turn to an important subject in probability and statistics: \emph{Copula}. A (bivariate) copula is a function $C:[0,1]^2\to[0,1]$ having the following properties:
\begin{enumerate} [$\qquad(1)$]
\item {$C(0,\cdot)\equiv0\equiv C(\cdot,0)$,}
\item {$C(u,1)=u$ and $C(1,v)=v$ for any $u,v\in[0,1]$, and}
\item {$C$ is quasi-monotone.}
\end{enumerate}
Copulas are a powerful tool in high-dimensional statistical applications as they have been widely used to describe the dependence relation between random variables. This is the result of the celebrated \emph{Sklar's Theorem}, stating that a bivariate distribution function $F$ is coupled with its one-dimensional marginals in the way that
\begin{equation} \label{copuladef}
F=C(F_1,F_2)
\end{equation}
for some copula $C$ \cite{Copula}. The same theorem also points out that $C$ for which (\ref{copuladef}) holds is unique if both the marginals of $F$ are continuous.

The tail distribution of $F$ is linked together with the marginals by the following formula
\begin{equation} \label{survivalC1}
\overline{F}=\widehat{C}(\overline{F_1},\overline{F_2}),
\end{equation}
where
\begin{equation} \label{survivalC2}
\widehat{C}(u,v)=C(\overline{u},\overline{v})+u+v-1
\end{equation}
is also a copula, called the \emph{survival copula} associated with $C$. Here the notational convention $\overline{a}=1-a$ for $a\in[0,1]$ is adopted. Note that a copula $C$ itself is the survival copula associated with $\widehat{C}$.

An application of Remark \ref{qmonotone} indicates that a copula is increasing and uniformly continuous on its domain:
\begin{equation} \label{Cunifconti}
0\leq C(u_2,v_2)-C(u_1,v_1)\leq(u_2-u_1)+(v_2-v_2),\;\;\;\;\;(u_1,v_1)\leq(u_2,v_2).
\end{equation}
We also have $C(u,v)\leq\min\{u,v\}$ on $[0,1]^2$, where $\min\{u,v\}$ is known as the \emph{comonotone copula}.

While bivariate copulas can be used to build bivariate distribution functions, the extension of univariate extreme value analysis to the bivariate framework requires a particular family of copulas, called extreme-value copulas. These copulas arise as limiting copulas of suitably normalized componentwise maxima of independent and identically distributed sequences of random vectors. There have been extensive works on extreme-value copulas in the literature, e.g. see the monograph \cite{Joe}. For readers' convenience and the sake of completeness, we list below some basics needed for the use in the paper.

A copula $C^*$ is an extreme-value copula if there exists a copula $C$ such that as $n\to\infty$,
\begin{equation} \label{extreCdef1}
C^n(u^{1/n},v^{1/n})\to C^*(u,v),\;\;\;\;\;(u,v)\in[0,1]^2,
\end{equation}
in which case $C$ is said to be in the \emph{domain of attraction} of $C^*$.
By applying a linear expansion, one can easily see that (\ref{extreCdef1}) is equivalent to
\begin{equation} \label{extreCdef2}
\lim_{\epsilon\to0^+}\epsilon^{-1}\big[1-C(1-\epsilon x_1,1-\epsilon x_2)\big]=-\log C^*(e^{-x_1},e^{-x_2})
\end{equation} for $(x_1,x_2)\in[0,\infty)^2$. An extreme-value copula is max-stable, i.e., the identity $C^*(u^{1/n},v^{1/n})^n=C^*(u,v)$ holds on $[0,1]^2$ for any $n\in\mathbb{N}$ by definition, and vice versa.

Analytically, a bivariate extreme-value copula is uniquely
characterized by means of a finite measure $\rho$ on the unit simplex $S_2=\{(x,y)\in[0,1]^2:x+y=1\}$, referred to as \emph{Pickands dependence measure}, which is a result due to De Haan and Resnick \cite{Haan} and Pickands \cite{Pickands}. More precisely, a copula is an extreme-value copula if and only if it is of the form
\begin{equation} \label{extremeC1}
C_A^*(u,v)=\exp\left[\log(uv)A\left(\frac{\log u}{\log(uv)}\right)\right],\;\;\;\;\;(u,v)\in(0,1)^2,
\end{equation}
where
\begin{equation} \label{Arepre}
A(t)=\int_{S_2}\max\{tx,(1-t)y\}\;d\rho(x,y)
\end{equation} and the aforementioned measure
satisfies the mean constraints
\begin{equation} \label{rhocond}
\int_{S_2}x\;d\rho(x,y)=1=\int_{S_2}y\;d\rho(x,y).
\end{equation}
This function $A:[0,1]\to[1/2,1]$, known as the \emph{Pickands dependence function}, is convex and satisfies the inequalities $t\vee(1-t)\leq A(t)\leq1$ on $[0,1]$.

The upper bound $A\equiv1$ corresponds to the independence copula $\Pi(u,v)=uv$, while the lower
bound $A(t)=t\vee(1-t)$ corresponds to the comonotone copula.

\subsection{Extremes and Max-Infinite Divisibility} Let $\mathbf{X}^{(i)}$, $i=1,\ldots,n$, be independent and identically distributed random vectors of $\mathbb{R}^2$ having common distribution function $H$. For $j=1,2$, let $M_j^{(n)}=\max\{X_j^{(1)},\ldots,X_j^{(n)}\}$. In extreme value analysis, one seeks for well chosen normalizing $\mathbb{R}^2$-sequences
$(\mathbf{a}^{(i)})_i$ and $(\mathbf{b}^{(i)})_i$ with $a_j^{(i)}>0$ and a proper bivariate distribution function $G$
with non-degenerate marginals so that as $n\to\infty$,
\begin{equation} \label{EV1}
\mathrm{Pr}\left[\frac{M_j^{(n)}-b_j^{(n)}}{a_j^{(n)}}\leq x_j\right]=H^n\big(a_1^{(n)}x_1+b_1^{(n)},a_2^{(n)}x_2+b_2^{(n)}\big)\to G(\mathbf{x})
\end{equation}
weakly. The limit (\ref{EV1}) is equivalently rephrased as
\begin{equation} \label{EV2}
\lim_{n\to\infty}n\big[1-H\big(a_1^{(n)}x_1+b_1^{(n)},a_2^{(n)}x_2+b_2^{(n)}\big)\big]=-\log G(\mathbf{x})
\end{equation} for any $\mathbf{x}\in\{G>0\}$ \cite{characterG}. Any $H$ giving rise to (\ref{EV1}) will be said to be in the max-domain of attraction of $G$ and this will be written $H\in\mathcal{D}_*(G)$.

Each marginal of the limiting distribution function $G$ is either the type of Gumbel $(\xi=0)$, Fr\'{e}chet $(\xi>0)$ or Weibull $(\xi<0)$:
\[G_j(x)=\exp\left[-\left(1+\xi\cdot\frac{x-m}{\sigma}\right)^{-1/\xi}\right],\;\;\;\;\;1+\xi(x-m)/\sigma>0,\] where $(\xi,m,\sigma)\in\mathbb{R}\times\mathbb{R}\times(0,\infty)$.
A bivariate extreme value distribution $G$, unlike the univariate case, cannot be represented as a distribution function indexed by a finite-dimensional parameter vector. However, there is an elegant representation formula for $G$ through its marginals and a unique Pickands dependence function $A$:
\begin{equation} \label{Gextreme}
G=\exp\left[\big(\log G_1+\log G_2\big)A\left(\frac{\log G_1}{\log G_1+\log G_2}\right)\right]
\end{equation}
whenever $G>0$ \cite{Grepre}.

A distribution function $G$ is max-stable if for $j=1,2$ and every $n\in\mathbb{N}$, there exist real numbers
$a_j^{(n)}>0$ and $b_j^{(n)}$ such that
\begin{equation} \label{maxstable}
G^n\big(a_1^{(n)}x_1+b_1^{(n)},a_2^{(n)}x_2+b_2^{(n)}\big)=G(\mathbf{x}).
\end{equation}
Clearly, a max-stable distribution function is in its own max-domain of attraction. As a matter of fact, the class of bivariate extreme value distributions is precisely
the class of max-stable distribution functions with non-degenerate marginals.

An immediate consequence of (\ref{maxstable}) is that $G^{1/n}$ is a distribution function for every positive integer $n$, meaning that $G$ is max-infinitely divisible according to the definition. It was shown in \cite{Resnick} that there exists a $\sigma$-measure $\tau$ on $[-\boldsymbol\infty,\boldsymbol\infty)$ such that
\begin{equation} \label{exponent}
G(\mathbf{x})=\exp\big\{-\tau\big([-\boldsymbol\infty,\boldsymbol\infty)
\backslash[-\boldsymbol\infty,\mathbf{x}]\big)\big\},\;\;\;\;\; \mathbf{x}\in[-\boldsymbol\infty,\boldsymbol\infty),
\end{equation}
whence the name exponent measure. Undoubtedly, any distribution function written as form (\ref{exponent}) is max-infinitely divisible.

\subsection{Free Extremes} We now briefly present the free extremal types theorem. Given a sequence of free self-adjoint operators $(X_n)$ having the common spectral probability distribution on the real line, one asks for suitably normalized constants $a_n>0$ and $b_n$ ensuring the weak convergence of
\begin{equation} \label{freemax}
\bigvee(X_n-b_n I)/a_n.
\end{equation}
Here, the notation $\vee$ stands for the max-operation on self-adjoint operators with respect to the spectral order.

One crucial fact underlying the asymptotic question in (\ref{freemax}) is that it can be converted to classical analysis problems by working with the free max-convolution of univariate distribution functions and with the map
\begin{equation} \label{freemap}
x\mapsto(1+\log x)_+,\;\;\;\;\;x\in\mathbb{R},
\end{equation} where
$c_+=\max\{c,0\}$ for $c\in\mathbb{R}$ \cite{freeextre}. The free max-convolution of two univariate distribution functions $F_1$ and $F_2$, denoted by $F_1\maxconv F_2$, can be computed through the concise formula
\[F_1\maxconv F_2=(F_1+F_2-1)_+.\] Furthermore, the map in (\ref{freemap}) serves as a homomorphism from the semigroup of univariate distribution functions endowed with the pointwise multiplication to the same set endowed with the free max-convolution. Consequently, one obtains the bijective correspondence between classical and free extremal types theorem via this homomorphism.

As in the classical case, a univariate distribution function is freely max-stable if and only if it is of a free extremal type, that is, it is the same type as one of the exponential distribution, Pareto distribution, and the Beta law. Another rigid tie between the classical extreme value theory and its free counterpart is that the classical and free max-domains of attraction of the corresponding laws not only coincide, but also share the identical normalizing constants \cite{freeextre}.

\subsection{Bi-free Max-Convolution} Bi-free max-convolution of compactly supported planar probability distributions was introduced in \cite{V15EXT} from the perspective of bounded linear operators. By working with affiliated unbounded self-adjoint operators, this convolution operation extends to planar probability measures without compact support. Below we adopt the perspective of distribution functions to quickly explain how to eliminate this compactness limitation based on the work in \cite{V15EXT}.

Given two arbitrary bivariate distribution functions $F$ and $G$, let $F^{(n)}$ and $G^{(n)}$ be distribution functions of compactly supported probability distributions on $\mathbb{R}^2$ converging weakly to $F$ and $G$, respectively. That is, $F^{(n)}\to F$ on $\mathcal{C}(F)$ and $G^{(n)}\to G$ on $\mathcal{C}(G)$ pointwise.

Denote by $H^{(n)}$ the bi-free max-convolution of $F^{(n)}$ and $G^{(n)}$ \cite{V15EXT}. In other words, $H^{(n)}$ is the distribution function of certain planar probability distribution $\mu^{(n)}$ with compact support, whose marginals are given by the formula
\begin{equation} \label{maxconvdef1}
H^{(n)}_j=\big(F^{(n)}_j+G^{(n)}_j-1\big)_+,\;\;\;\;\;j=1,2,
\end{equation}
while $H^{(n)}(\mathbf{x})$ itself is determined by the identity
\begin{equation} \label{maxconvdef2}
\frac{H^{(n)}_1H^{(n)}_2}{H^{(n)}}=\frac{F^{(n)}_1F^{(n)}_2}{F^{(n)}}+
\frac{G^{(n)}_1G^{(n)}_2}{G^{(n)}}-1
\end{equation} provided that none of $F^{(n)}(\mathbf{x})$, $G^{(n)}(\mathbf{x})$, $H_1^{(n)}(x_1)$, and $H_2^{(n)}(x_2)$ is zero and $H^{(n)}(\mathbf{x})=0$ elsewhere.

Each family $\{\mu_j^{(n)}\}$ of marginals is tight in view of (\ref{maxconvdef1}), whence so is the family $\{\mu^{(n)}\}$ itself. Let $H$ be a weak-limit distribution function of $H^{(n)}$. For any point $\mathbf{x}$ with $x_j\in\mathcal{C}(F_j)\cap\mathcal{C}(G_j)\cap\mathcal{C}(H_j)$ for $j=1,2$ (such an $\mathbf{x}$ also lies in $\mathcal{C}(F)\cap\mathcal{C}(G)\cap\mathcal{C}(H)$ according to (\ref{upperbound})),
letting $n\to\infty$ in (\ref{maxconvdef1}) gives that
\begin{equation} \label{Hmarginal}
H_j=\big(F_j+G_j-1\big)_+
\end{equation} holds at $x_j$. If, in addition, $\mathbf{x}$ is also selected from the positive set $\mathcal{P}=\{F>0\}\cap\{G>0\}\cap\{H_1>0\}\times\{H_2>0\}$, then
letting $n\to\infty$ in (\ref{maxconvdef2}) results that $H(\mathbf{x})>0$ and the following identity is valid at $\mathbf{x}$:
\begin{equation} \label{H}
\frac{H_1H_2}{H}=\frac{F_1F_2}{F}+
\frac{G_1G_2}{G}-1.
\end{equation}

Since the discontinuities of an increasing function on the real line are at most countable, the right-continuity of univariate distribution functions implies that any weak-limit of $H^{(n)}$ has the same marginals, which are clearly given by the equation (\ref{Hmarginal}) on $\mathbb{R}$. The above-continuity of bivariate distribution functions then yields that $H$ is uniquely determined by the equation (\ref{H}) on the set $\mathcal{P}$. Apparently, both equations (\ref{Hmarginal}) and (\ref{H}) are independent of the choice of the weak-convergent distribution functions $F^{(n)}$ and $G^{(n)}$. These findings consequently lead to the following definition.

\begin{pdef} \emph{The \emph{bi-free max-convolution} $H$ of two bivariate distribution functions $F$ and $G$, denoted by
\[F\maxconv\maxconv G,\] is the unique bivariate distribution function so that its marginals satisfy the equation (\ref{Hmarginal}) on $\mathbb{R}$ for $j=1,2$, and $H(\mathbf{x})$ itself is determined by (\ref{H}) whenever $F(\mathbf{x})$, $G(\mathbf{x})$, $H_1(x_1)$, and $H_2(x_2)$ are all strictly positive and $H=0$ elsewhere. The
bi-free max-convolution $\mu\maxconv\maxconv\nu$ of two planar probability distributions $\mu$ and $\nu$ having distribution functions $F$ and $G$, respectively, is the planar probability distribution having $F\maxconv\maxconv G$ as its distribution function.}
\end{pdef}

We remark here that the marginals of $F\maxconv\maxconv G$ are the free max-convolution of $F_j$ and $G_j$:
\[(F\maxconv\maxconv G)_j=F_j\maxconv G_j,\] which is defined as in (\ref{Hmarginal}), e.g. see \cite{freeextre}.
One can also easily verify from the definition the associativity of the binary operation $\maxconv\maxconv$ for arbitrary distribution functions $F$, $G$ and $H$:
\[(F\maxconv\maxconv\;G)\maxconv\maxconv\;H=F\maxconv\maxconv(G\maxconv\maxconv\;H).\]

\section{Limit Theorems and Bi-Freely max-infinite divisibility} Having defined the bi-free max-convolution of bivariate distribution functions, we can now introduce the first main topic under the study in the paper.

\begin{pdef} \emph{A distribution function $F$ on $\mathbb{R}^2$ is said to be \emph{bi-freely max-infinitely divisible} (abbreviated to \emph{bi-freely max-i.d.}) if for any $n\in\mathbb{N}$, there exists one distribution function $F^{(1/n)}$ on $\mathbb{R}^2$ so that
\begin{equation} \label{bifreemax}
(F^{(1/n)})^{\maxconvexp\maxconvexp n}:=\underbrace{F^{(1/n)}\maxconv\maxconv\cdots\maxconv\maxconv\; F^{(1/n)}}_{n\;\;\mathrm{terms}}=F.
\end{equation} The planar probability distribution $\mu^{\maxconvexp\maxconvexp n}$ and the bi-freely max-infinite divisibility of a planar probability distribution $\mu$ are defined analogously. Denote by $\mathcal{ID}(\maxconv\maxconv)$ the class of bi-freely max-i.d. distribution functions (or bi-freely max-i.d. probability measures).}
\end{pdef}

For ease of reference, we restate (\ref{bifreemax}) as follows: $F_j=(nF_j^{(1/n)}-(n-1))_+$ holds on $\mathbb{R}$ for $j=1,2$, and $F(\mathbf{x})$ satisfies the identity
\begin{equation} \label{maxFexpress}
\frac{F_1F_2}{F}-1=n\left[\frac{F_1^{(1/n)}F_2^{(1/n)}}{F^{(1/n)}}-1\right]
\end{equation}
if none of $F^{(1/n)}(\mathbf{x})$, $F_1(x_1)$ and $F_2(x_2)$ vanishes and $F(\mathbf{x})=0$ otherwise.

The results in the following remark will be frequently employed in the paper.

\begin{remark} \label{maxremark}
\emph{The necessary and sufficient conditions for (\ref{bifreemax}) to be true for some $F^{(1/n)}$ are:
\begin{enumerate} [$\qquad(1)$]
\item {$F_j=nF^{(1/n)}_j-(n-1)$ holds on the positive set $\mathcal{P}_j=\{F_j>0\}$ for $j=1,2$,}
\item {$\{F>0\}=\mathcal{P}_1\times\mathcal{P}_2$, and}
\item {$F$ satisfies the identity (\ref{maxFexpress}) whenever it is nonzero.}
\end{enumerate} The necessity follows by observing that we have $F^{(1/n)}>1-2/n$ on $\mathcal{P}_1\times\mathcal{P}_2$, which is due to the inequalities
$F^{(1/n)}+\alpha+\beta\leq1$, $F^{(1/n)}+\alpha>1-1/n$ and $F^{(1/n)}+\beta>1-1/n$ on $\mathcal{P}_1\times\mathcal{P}_2$, where $\alpha=F^{(1/n)}_1-F^{(1/n)}$ and $\beta=F^{(1/n)}_2-F^{(1/n)}$. If conditions (1)-(3) hold, then redefining $F^{(1/n)}$ on $\mathbb{R}\backslash\mathcal{P}_1\times\mathcal{P}_2$, redefining as zero for instance, so that the statement in (1) is still valid confirms the sufficiency. These discussions also come to that no unique $F^{(1/n)}$ satisfies (\ref{bifreemax}).}
\end{remark}

The following result is a direct consequence of the associativity and commutativity of the binary operation $\maxconv\maxconv$ on bivariate distribution functions.

\begin{prop} \label{closed}
The bi-free max-convolution of two bi-freely max-infinitely divisible distribution functions is again bi-freely max-infinitely divisible.
\end{prop}

Paralleling the classical theory, the set $\mathcal{ID}(\maxconv\maxconv)$ is closed under the topology of weak convergence, which is obtained by the following theorem.

\begin{thm} \label{maxconv1}
Let $F,F^{(1)},F^{(2)},\ldots$ be bivariate distribution functions and let
$\mathcal{C}_j=\mathcal{C}(F_j)\cap\{F_j>0\}$ for $j=1,2$. If $\inf\{F_j>0\}>-\infty$ for $j=1,2$, and if as $n\to\infty$,
\[\big(F_j^{(n)}\big)^{\maxconvexp\;n}\to F_j\;\;\;on\;\;\;\mathcal{C}_j\] and
\[\big(F^{(n)}\big)^{\maxconvexp\maxconvexp\;n}\to F\;\;\;on\;\;\;\mathcal{C}_1\times\mathcal{C}_2,\] then $F\in\mathcal{ID}(\maxconv\maxconv)$.
\end{thm}

\begin{pf} Observe first that the conditions listed in the proposition are equivalent to that as $n\to\infty$,
\begin{equation} \label{maxconv1eq1}
\big(nF^{(n)}_j-(n-1)\big)_+\to F_j\;\;\;on\;\;\;\mathcal{C}_j
\end{equation}
for $j=1,2$, and
\begin{equation} \label{maxconv1eq2}
n\left[\frac{F^{(n)}_1F^{(n)}_2}{F^{(n)}}-1\right]\to\frac{F_1F_2}{F}-1
\;\;\;on\;\;\;\mathcal{C}_1\times\mathcal{C}_2.
\end{equation}
Hence we may assume that each $F_j^{(n)}$ also vanishes on $\{F_j=0\}$ because
the conditions (\ref{maxconv1eq1}) and (\ref{maxconv1eq2}) are not affected at all by doing so. For convention, the term $F_1F_2/F$ in (\ref{maxconv1eq2}) is realized as infinity at points in $\mathcal{C}_1\times\mathcal{C}_2$ at which $F$ vanishes. We shall prove below, in fact, that the relation $\{F>0\}=\{F_1>0\}\times\{F_2>0\}$ is valid provided that (\ref{maxconv1eq1}) and (\ref{maxconv1eq2}) are satisfied.

Let $k\geq2$ be an arbitrary but fixed integer. Further let $G^{(n)}=(F^{(nk)})^{\maxconvexp\maxconvexp\;n}$, i.e., $G^{(n)}$ is the bivariate distribution function so that
$G^{(n)}_j=(nF^{(nk)}_j-(n-1))_+$ for $j=1,2$, and $G^{(n)}$ satisfies
\begin{equation} \label{limitcond1}
\frac{G^{(n)}_1G^{(n)}_2}{G^{(n)}}-1=n\left[\frac{F^{(nk)}_1F^{(nk)}_2}{F^{(nk)}}-1\right]
\end{equation}
on the set $\{F^{(nk)}>0\}\cap\{G_1^{(n)}>0\}\times\{G_2^{(n)}>0\}$ and $G^{(n)}=0$ elsewhere.

One can infer from (\ref{maxconv1eq1}) that
\begin{align*}
\big(kG^{(n)}_j-(k-1)\big)_+&=\max\big\{k\cdot\max\big\{nF_j^{(nk)}-(n-1),0\big\}-(k-1),0\big\} \\
&=\max\big\{nkF_j^{(nk)}-(nk-1),-k+1,0\big\} \\
&=\big(nkF^{(nk)}_j-(nk-1)\big)_+ \\
&\to F_j.
\end{align*} on $\mathcal{C}_j$ as $n\to\infty$. This together with the hypothesis $\inf\{F_j^{(nk)}>0\}=\inf\{F_j>0\}>-\infty$ ensures the tightness of the probability measures having $G_j^{(n)}$ as distribution functions. By means of (\ref{maxconv1eq2}) and (\ref{limitcond1}), it follows that
\begin{equation} \label{limitcond2}
k\left[\frac{G^{(n)}_1G^{(n)}_2}{G^{(n)}}-1\right]\to\frac{F_1F_2}{F}-1
\end{equation}
on $\mathcal{C}_1\times\mathcal{C}_2$ as $n\to\infty$ because of $\mathcal{C}_j=\mathcal{C}(F_j^{(nk)})\cap\{F_j^{(nk)}>1-1/(nk)\}=\mathcal{C}(G_j^{(n)})\cap\{G_j^{(n)}>1-1/k\}$ for all $n$.

Next, let $G$ be a weak-limit distribution function of the family $\{G^{(n)}\}$. Then by passing to a subsequence of $\{G^{(n)}\}$ if needed, we see that $(kG^{(n)}_j-(k-1))_+\to(kG_j-(k-1))_+$ on $\mathcal{C}(G_j)$ as $n\to\infty$ for $j=1,2$, and
\begin{equation} \label{limitcond3}
k\left[\frac{G^{(n)}_1G^{(n)}_2}{G^{(n)}}-1\right]\to k\left[\frac{G_1G_2}{G}-1\right]
\end{equation}
on $\mathcal{C}(G_1)\times\mathcal{C}(G_2)\cap\{G>0\}$ as $n\to\infty$. Consequently, we come to the relation
\begin{equation} \label{limitcond4}
F_j=kG_j-(k-1)
\end{equation}
on $\mathcal{P}_j=\{F_j>0\}$. Also observe that (\ref{limitcond4}) leads to $G_j>1-1/k$ on $\mathcal{P}_j$, and so we derive the inclusion $\mathcal{P}_1\times\mathcal{P}_2\subset\{G>0\}$. These observations, along with (\ref{limitcond2}) and (\ref{limitcond3}), yield that
\begin{equation} \label{limitcond5}
k\left[\frac{G_1G_2}{G}-1\right]=\frac{F_1F_2}{F}-1
\end{equation} holds on the set $\mathcal{C}_1\times\mathcal{C}_2$, and on the set $\mathcal{P}_1\times\mathcal{P}_2$ as well by employing the (above-) right-continuity of distribution functions.

Finally, the validity of the formula (\ref{limitcond5}) on $\mathcal{P}_1\times\mathcal{P}_2$ shows that $\{F>0\}=\mathcal{P}_1\times\mathcal{P}_2$. Since the same reasonings also apply to the other $k\geq2$, we have verified the bi-freely max-infinite divisibility of $F$ by virtue of Remark \ref{maxremark}.
\end{pf} \qed

To proceed further, two fundamental and useful lemmas are needed.

\begin{lem} \label{ineq}
Let $\alpha,\beta,\gamma$, and $\delta$ be real numbers so that $\alpha\leq\delta$ and $\beta,\gamma\in[\alpha,\delta]$.
\begin{enumerate} [$\qquad(1)$]
\item {If $\delta-\gamma-\beta+\alpha\geq0$ \emph{(}resp. $>0$\emph{)}, then
$e^\delta-e^\gamma-e^\beta+e^\alpha\geq0$ \emph{(}resp. $>0$\emph{)}.}
\item {If $\alpha>0$ and $\delta-\gamma-\beta+\alpha\leq0$ \emph{(}resp. $<0$\emph{)}, then
$\log\delta-\log\gamma-\log\beta+\log\alpha\leq0$ \emph{(}resp. $<0$\emph{)}.}
\end{enumerate}
\end{lem}

\begin{pf} The results in (1) are an immediate consequence of the decomposition
\[e^\delta-e^\gamma-e^\beta+e^\alpha=e^\alpha(e^{\beta-\alpha}-1)(e^{\gamma-\alpha}-1)
+e^{\beta+\gamma-\alpha}
(e^{\delta-\gamma-\beta+\alpha}-1),\] while the conclusions in (2) are a direct application of (1).
\end{pf} \qed

\begin{lem} \label{tauform}
Let $\mathbf{L}\in\mathbb{R}^2$ and $\tau$ be a positive Borel measure on $\mathbb{R}^2$ with marginals $\tau_1$ and $\tau_2$. Suppose that for $j=1,2$, the function
\begin{equation} \label{tauformFj}
F_j(x)=\left\{
\begin{array}{ll}
1-\tau_j((x,\infty)), & \hbox{$x\geq L_j$,} \\
0, & \hbox{$x<L_j$,}
\end{array}
\right.
\end{equation} is a univariate distribution function having support $[L_j,\infty)$.
Then the function
\begin{equation} \label{tauformF}
F(\mathbf{x})=\frac{F_1(x_1)F_2(x_2)}{1-\tau((\mathbf{x},\boldsymbol\infty))},\;\;\;\;\;
\mathbf{x}>\mathbf{L},
\end{equation}
extends as a bivariate distribution function with $F_1$ and $F_2$ as its marginals and satisfies
$\{F>0\}=\{F_1>0\}\times\{F_2>0\}$. Moreover, $\log F$ is quasi-monotone on $(\mathbf{L},\boldsymbol\infty)$.
\end{lem}

\begin{pf} For notational simplicity, we denote $g(\mathbf{x})=1-\tau((\mathbf{x},\boldsymbol\infty))$ for $\mathbf{x}>\mathbf{L}$. Notice that on account of $0<F_1F_2\leq g$, $F$ is well-defined and $0<F\leq1$ on $(\mathbf{L},\boldsymbol\infty)$.

We first prove that $F(\mathbf{x})\leq F(\mathbf{y})$ for any points $\mathbf{x}\leq\mathbf{y}$ in $(\mathbf{L},\boldsymbol\infty)$. It is enough to consider the case $\mathbf{y}=(y_1,x_2)$ with $y_1>x_1$ as an analogous inequality holds for $\mathbf{y}=(x_1,y_2)$ with $y_2>x_2$. Since the difference
$F_1(y_1)g(\mathbf{x})-F_1(x_1)g(\mathbf{y})$ simplifies into
\[g(\mathbf{x})\tau((x_1,y_1]\times\mathbb{R})-\big[1-\tau((x_1,\infty)\times\mathbb{R})\big]
\tau((x_1,y_1]\times(x_2,\infty)),\] which is apparently non-negative, we gain the inequality $F(\mathbf{x})\leq F(\mathbf{y})$.

The verification of the quasi-monotonicity of $F$ and $\log F$ on $(\mathbf{L},\boldsymbol\infty)$ carries out as follows. Since we have $g(\mathbf{x})\leq g(x_1,y_2),g(y_1,x_2)\leq g(\mathbf{y})$ and $V_g([\mathbf{x},\mathbf{y}])=-\tau((x_1,y_1]\times(x_2,y_2])\leq0$ for points $\mathbf{x}\leq\mathbf{y}$ in $(\mathbf{L},\boldsymbol\infty)$, Lemma \ref{ineq} says that
\[V_{\log F}([\mathbf{x},\mathbf{y}])=-\log g(\mathbf{y})+\log g(x_1,y_2)+\log g(y_1,x_2)-\log g(\mathbf{x})\geq0\] and
\begin{equation} \label{VF}
V_F([\mathbf{x},\mathbf{y}])=e^{\log F(\mathbf{y})}-e^{\log F(x_1,y_2)}-e^{\log F(y_1,x_2)}+e^{\log F(\mathbf{x})}\geq0.
\end{equation}
Then the quasi-monotonicity of $F$ and the inequality (\ref{quasiineq}) enable us to derive
\[|F(\mathbf{y})-F(\mathbf{x})|\leq|F_1(y_1)-F_1(x_1)|+|F_2(y_2)-F_2(x_2)|,
\;\;\;\;\;\mathbf{x},\mathbf{y}\in(\mathbf{L},\boldsymbol\infty).\]

The inequality above clearly yields that $F$ has an extension to $[\mathbf{L},\boldsymbol\infty)$, and the extension is increasing, continuous from above, and quasi-monotone. By assigning $F$ the value $0$ on the complement $\mathbb{R}^2\backslash[\mathbf{L},\boldsymbol\infty)$, one can see that $F$ is a bivariate distribution function with the stated properties in the lemma.
\end{pf} \qed

We conclude this section with an example, which views as an analog of bi-free compound Poisson distribution \cite{HW1}. It will be shown in the next section that virtually all bi-freely max-i.d. distribution functions appear in the manner given in this example.

\begin{exam} \label{Poissonexam}
\emph{(Bi-free Max Analogue of Compound Poisson distribution). Given a number $\lambda>0$, a point $\mathbf{p}\in\mathbb{R}^2$, and a planar probability distribution $\nu$, define $\tau=\lambda\nu$ and let $F^{(n)}$ be the distribution function of the planar probability law
\[\left(1-\frac{\lambda}{n}\right)\delta_\mathbf{p}+\frac{\lambda}{n}\nu\] and $L_j=\max\big\{p_j,\inf\{F_{\nu_j}>1-1/\lambda\}\big\}$. As recorded in Lemma \ref{tauform}, one can use $\tau$ to construct a bivariate distribution function $F$ as in (\ref{tauformF}), where each marginal $F_j$ has support $[L_j,\infty)$ and has tail given by (\ref{tauformFj}). Then according to the construction, it is fairly easy to see that $(nF_j^{(n)}-(n-1))_+\to F_j$ on $\mathbb{R}$ as $n\to\infty$ for $j=1,2$. Furthermore, for any point $\mathbf{x}\in[\mathbf{L},\boldsymbol\infty)$, we have as $n\to\infty$,
\[n\left[\frac{F_1^{(n)}F_2^{(n)}}{F^{(n)}}-1\right](\mathbf{x})
=\lambda\big(F_{\nu_1}+F_{\nu_2}-F_\nu-1\big)(\mathbf{x})+O(n^{-1})\to
-\tau((\mathbf{x},\boldsymbol\infty)),\] which coincides with $(F_1F_2/F)(\mathbf{x})-1$. As a consequence of Theorem \ref{maxconv1}, $F$ is bi-freely max-i.d. and
\[\left(\left(1-\frac{\lambda}{n}\right)\delta_\mathbf{p}
+\frac{\lambda}{n}\nu\right)^{\maxconvexp\maxconvexp\;n}
\Rightarrow\mu,\] where $\mu$ is the planar probability measure having distribution function $F$.}
\end{exam}

\section{Characterization of bi-freely max-i.d. distribution functions} A distribution function $F$ on $\mathbb{R}^2$ is said to have support bounded from below if $\{F>0\}\subset[\mathbf{L},\boldsymbol\infty)$ for some $\mathbf{L}\in\mathbb{R}^2$. Further speaking, the positive set $\{F>0\}$ is a generalized rectangle of form $(\mathbf{L},\boldsymbol\infty)$, $(L_1,\infty)\times[L_2,\infty)$, $[L_1,\infty)\times(L_2,\infty)$ or $[\mathbf{L},\boldsymbol\infty)$ with $\mathbf{L}\in\mathbb{R}^2$ if and only if $F$ has support bounded from below and satisfies the relation
\begin{equation} \label{F>0}
\{F>0\}=\{F_1>0\}\times\{F_2>0\}.
\end{equation} This could be seen from the fact that $\mathbf{x},\mathbf{y}\in\{F>0\}$ implies $F(\mathbf{x}\wedge\mathbf{y})>0$ if $\{F>0\}$ is a generalized rectangle. We will denote the point $\mathbf{L}$ by $\inf\mathrm{supp}(F)$ when $F$ has support bounded from below and (\ref{F>0}) holds.

In this section, we shall characterize the collection $\mathcal{ID}(\maxconv\maxconv)$ of bi-freely max-i.d. distribution functions. We first introduce a transformation of the bivariate distribution function $F$:
\begin{equation} \label{TF}
T_F=\frac{F_1F_2}{F}-F_1-F_2+1
\end{equation} defined on the set $\{F>0\}$.
Clearly, $T_F\geq0$ on $\{F>0\}$. More information concerned with $T_F$ is described below.

\begin{lem} \label{main1}
The followings hold for an $F\in\mathcal{ID}(\maxconv\maxconv)$.
\begin{enumerate} [$\qquad(1)$]
\item {$F$ has support bounded from below and satisfies $\{F>0\}=\{F_1>0\}\times\{F_2>0\}$.}
\item {$-T_F$ is increasing, continuous from above,
and quasi-monotone on $\{F>0\}$.}
\end{enumerate}
\end{lem}

\begin{pf} Let $\Omega=\{F>0\}$ and let $F^{(1/n)}$ be a distribution function so that $(F^{(1/n)})^{\maxconvexp\maxconvexp n}=F$ for all $n$. That we have $\Omega=\{F_1>0\}\times\{F_2>0\}$ was already indicated in Remark \ref{maxremark}. If the support of $F$ is not bounded from below, for instance, if $F_1$ never vanishes on $\mathbb{R}$, then again in view of Remark \ref{maxremark}, we will have $F_1^{(1/2)}>1/2$ on $\mathbb{R}$, a contradiction. Hence (1) is established.

To prove the properties regarding the function $T_F$ in assertion (2), we first express
\[F^{(1/n)}=\frac{(F_1+n-1)(F_2+n-1)}{n(T_F+F_1+F_2+n-2)}\]
on $\Omega$. This expression leads to that for any points $\mathbf{x}\leq\mathbf{y}$ in $\Omega$, we have
\begin{equation} \label{TFdiffer1}
0\leq n\big[F^{(1/n)}(\mathbf{y})-F^{(1/n)}(\mathbf{x})\big]=T_F(\mathbf{x})-T_F(\mathbf{y})+O(n^{-1})
\end{equation}
as $n\to\infty$. Since $F_j=(nF_j^{(1/n)}-(n-1))_+$, the use of (\ref{TFdiffer1}) and (\ref{upperbound}) then reveals that as $n\to\infty$,
\begin{align*}
0&\leq T_F(\mathbf{x})-T_F(\mathbf{y})+O(n^{-1}) \\
&\leq n\left[F_1^{(1/n)}(y_1)-F_1^{(1/n)}(x_1)+F_2^{(1/n)}(y_2)-F_2^{(1/n)}(x_2)\right] \\
&=F_1(y_1)-F_1(x_1)+F_2(y_2)-F_2(x_2).
\end{align*} This concludes that the inequalities
\begin{equation} \label{TFdiffer2}
0\leq T_F(\mathbf{x})-T_F(\mathbf{y})\leq F_1(y_1)-F_1(x_1)+F_2(y_2)-F_2(x_2),
\end{equation} are true for any $\mathbf{x}\leq\mathbf{y}\in\Omega$, proving that $-T_F$ increases and is continuous from above on $\Omega$.

In order to obtain the quasi-monotonicity of $-T_F$ on $\Omega$, one can replace $F^{(1/n)}(\mathbf{y})-F^{(1/n)}(\mathbf{x})$ in (\ref{TFdiffer1}) with the differences  $F^{(1/n)}(\mathbf{y})-F^{(1/n)}((x_1,y_2))$ and $F^{(1/n)}((y_1,x_2))-F^{(1/n)}(\mathbf{x})$ and apply the same technique as shown above. This finishes the proof of assertion (2).
\end{pf} \qed

\begin{remark} \label{QFremark}
\emph{Given an $F\in\mathcal{ID}(\maxconv\maxconv)$, one can employ the transform
\begin{equation} \label{QF}
Q_F(\mathbf{x})=\frac{F_1F_2}{F}(\mathbf{x}),\;\;\;\;\;\mathbf{x}\in\{F>0\},
\end{equation} instead of working with $T_F$, and obtain its properties equivalent to those of $T_F$ in Lemma \ref{main1}. Specifically speaking, the function $-Q_F$ is quasi-monotone on $\{F>0\}$ and satisfies
\begin{equation} \label{QFprop}
0\leq Q_F(\mathbf{y})-Q_F(\mathbf{x})\leq F_1(y_1)-F_1(x_1)+F_2(y_2)-F_2(x_2)
\end{equation}
for any points $\mathbf{x}\leq\mathbf{y}$. Therefore, one may investigate $F$ via $T_F$ or $Q_F$, depending on which is more convenient.}
\end{remark}

Note that the distribution function $F$ of an almost surely constant random vector $\mathbf{X}$ in $\mathbb{R}^2$ is automatically bi-freely max-i.d. due to $F^{\maxconvexp\maxconvexp n}=F$ for all $n\in\mathbb{N}$. Such an $F$ satisfies $T_F=0$ on $[\mathbf{L},\boldsymbol\infty)$ if $\mathbf{X}=\mathbf{L}$ a.s. In the sequel, we shall concentrate on the distribution functions of random vectors that are not a.s. constant.

After the previous preparations, we are now ready to characterize the class $\mathcal{ID}(\maxconv\maxconv)$.

\begin{thm} \label{mainthm}
Suppose that $F$ is a bivariate distribution function, which is not the distribution function of an almost surely constant random vector, and that the set $\{F>0\}$ is a generalized rectangle with $\mathbf{L}=\inf\mathrm{supp}(F)>-\boldsymbol\infty$.
Then the following are equivalent.
\begin{enumerate} [$\qquad(1)$]
\item {The distribution function $F$ is bi-freely max-infinitely divisible.}
\item {The function $-Q_F$ is quasi-monotone on $\{F>0\}$
and satisfies \emph{(\ref{QFprop})} for any $\mathbf{x}\leq\mathbf{y}$.}
\item {The function $-T_F$ is increasing, continuous from above, and quasi-monotone on
$\{F>0\}$.}
\item {There exist a distribution function $H$ on $\mathbb{R}^2$ and
a finite number $\lambda>0$ so that
\begin{equation} \label{TFrepre1}
T_F=\lambda(1-H)
\end{equation}
holds on $\{F>0\}$.}
\item {There exists a positive Borel measure $\tau\neq0$ on $\mathbb{R}^2$ so that
\begin{equation} \label{Fjrepre}
\overline{F_j}(x)=\tau_j((x,\infty)),\;\;\;\;\;x>L_j,
\end{equation} and
\begin{equation} \label{Frepre}
F(\mathbf{x})=\frac{F_1(x_1)F_2(x_2)}{1-\tau((\mathbf{x},\boldsymbol\infty))},
\;\;\;\;\;\mathbf{x}>\mathbf{L}.
\end{equation}}
\end{enumerate}
Let $\nu$ be the probability measure having distribution function $H$. If \emph{(1)-(5)} hold, then $\tau=\lambda\nu$ and $T_F$ has an extension to $[\mathbf{L},\boldsymbol\infty)$ so that \emph{(\ref{TFrepre1})} holds on $[\mathbf{L},\boldsymbol\infty)$ and $0<T_F(\mathbf{L})\leq\lambda$.
\end{thm}

\begin{pf} (1)$\Rightarrow$(2) and (2)$\Rightarrow$(3): These two implications have been confirmed in Lemma \ref{main1} and Remark \ref{QFremark}.

(3)$\Rightarrow$(4): Assume that (3) holds. The inequality
\[|T_F(\mathbf{x})-T_F(\mathbf{y})|\leq|F_1(x_1)-F_1(y_1)|+|F_2(x_2)-F_2(y_2)|\] can be proved to be true for any $\mathbf{x},\mathbf{y}\in\{F>0\}$ directly (or using the equivalent of (2) and (3)). This estimate says that $T_F$ has an extension to $[\mathbf{L},\boldsymbol\infty)$ (the extension is also denoted by $T_F$), and the stated properties in (3) remain true for $T_F$ on $[\mathbf{L},\boldsymbol\infty)$ and $T_F(\mathbf{L})<\infty$. One also has $T_F(\mathbf{L})>0$, otherwise, $T_F(\mathbf{L})=0$ will imply $T_F\equiv0$, i.e., $F$ is the distribution function of the Dirac measure $\delta_\mathbf{L}$ concentrated at point $\mathbf{L}$. Now pick any number $\lambda\geq T_F(\mathbf{L})$. Then any bivariate distribution function consistent with $1-T_F/\lambda$ on $[\mathbf{L},\boldsymbol\infty)$ serves as an $H$ in (\ref{TFrepre1}). Hence (3)$\Rightarrow$(4) is established.

(4)$\Rightarrow$(5): Let $\nu$ be the probability measure having the distribution function $H$ and $\tau=\lambda\nu$. Then (\ref{TFrepre1}) shows that for any $x>L_1$,
\[\overline{F_1}(x)=\lim_{x_2\to\infty}T_F(x,x_2)=\lambda\overline{H_1}(x)=\tau_1((x,\infty)),\] which gives the representation (\ref{Fjrepre}) when $j=1$. Similarly, (\ref{Fjrepre}) is true for $j=2$. Note that (\ref{Fjrepre}) particularly shows that $\tau((\mathbf{x},\boldsymbol\infty))<1$ for $\mathbf{x}\in(\mathbf{L},\boldsymbol\infty)$. Since the identity $(1+H-H_1-H_2)(\cdot)=\nu((\cdot,\boldsymbol\infty))$ is valid on $(\mathbf{L},\boldsymbol\infty)$, (\ref{Frepre}) is obtained from (\ref{TFrepre1}), and so (4) implies (5).

(5)$\Rightarrow$(1): Assume that the statements in (5) are satisfied. To attain the bi-freely max-infinite divisibility of $F$, let $F^{(1/n)}$, $n\geq2$, be the bivariate distribution function with the properties that $\{F^{(1/n)}>0\}=[\mathbf{L},\boldsymbol\infty)$,
\[F_j^{(1/n)}(x)=1-n^{-1}\tau_j((x,\infty)),\;\;\;\;\;x>L_j,\] for $j=1,2$, and
\[F^{(1/n)}(\mathbf{x})=\frac{F^{(1/n)}_1(x_1)F^{(1/n)}_2(x_2)}
{1-n^{-1}\tau((\mathbf{x},\infty))},\;\;\;\;\;\mathbf{x}>\mathbf{L}.\] The existence of such an $F^{(1/n)}$ is guaranteed by Lemma \ref{tauform}. It is fairly easy to verify that $F_j=(nF_j^{(1/n)}-(n-1))_+$ is valid for $j=1,2$ and that $nT_{F^{(1/n)}}=T_F$ holds on the set $\{F>0\}$. Hence we come to the conclusion $(F^{(1/n)})^{\maxconvexp\maxconvexp n}=F$, and so (5)$\Rightarrow$(1) is proved.
\end{pf} \qed

The proof of Theorem \ref{mainthm} provides us with a partial semigroup $\{F^{(t)}\}_{t\geq0}$ of bivariate distribution functions generated by a bi-freely max-i.d. distribution function $F$:
\begin{equation} \label{semigroup}
T_{F^{(t)}}=tT_F
\end{equation} on $\{F>0\}\cap\{F^{(t)}>0\}$. To be precise, let $F$ be given by (\ref{Frepre}). Then for $t\in[0,1]$, any bivariate distribution function $F^{(t)}$ satisfying
\begin{equation} \label{semigroup1}
F_j^{(t)}(x_j)=tF_j(x_j)-(t-1),\;\;\;\;\;x_j\in\{F_j>0\},
\end{equation}
and
\begin{equation} \label{semigroup2}
F^{(t)}(\mathbf{x})=\frac{F^{(t)}_1(x_1)F^{(t)}_2(x_2)}
{1-t\tau((\mathbf{x},\boldsymbol\infty))},\;\;\;\;\;\mathbf{x}\in\{F>0\},
\end{equation} will do the job. The construction clearly shows that we have the inclusion $\mathrm{supp}(F)\subset\mathrm{supp}(F^{(t)})$. Unless the relation $\mathrm{supp}(F^{(t)})=\mathrm{supp}(F)$ is a requirement, when $t\in[0,1)$, $F^{(t)}$ satisfying the equation (\ref{semigroup}) in general is not unique.

As for the case of $t>1$, let $\mathbf{L}^{(t)}=(L^{(t)}_1,L_2^{(t)})$, where
$L_j^{(t)}=\inf\{F_j>1-1/t\}$. Then
the distribution function $F^{(t)}$ having support $[\mathbf{L}^{(t)},\boldsymbol\infty)$ and satisfying
\[F^{(t)}_j=(tF_j-(t-1))_+\] and
\[F^{(t)}(\mathbf{x})=\frac{F_1^{(t)}(x_1)F_2^{(t)}(x_2)}{1-t\tau((\mathbf{x},\boldsymbol\infty))},
\;\;\;\;\;\mathbf{x}\in(\mathbf{L}^{(t)},\boldsymbol\infty),\] is the desired one. Different from the previous situation, such an $F^{(t)}$ is unique and $\mathrm{supp}(F^{(t)})\subset\mathrm{supp}(F)$.

We close this section with some special classes of bi-freely max-i.d. distribution functions.

\begin{cor} \label{prod}
Let $F$ be a bivariate distribution function.
\begin{enumerate} [$\qquad(1)$]
\item {Suppose that $F\in\mathcal{ID}(\maxconv\maxconv)$ and $\inf\mathrm{supp}(F)=\mathbf{L}$.
Then $T_F(\mathbf{L})=2-F_1(L_1)-F_2(L_2)$ if and only if $F$
is the distribution function of two independent random variables, in which case $T_F=2-F_1-F_2$ on $[\mathbf{L},\boldsymbol\infty)$.}
\item {Suppose that $F$ is the distribution of two independent random variables.
Then $F\in\mathcal{ID}(\maxconv\maxconv)$ if and only
if its support is bounded from below.}
\end{enumerate}
\end{cor}

\begin{pf} Assertion (1) is apparently established if $F$ is the distribution function of $\delta_\mathbf{L}$. For other $F$, by virtue of the representation (\ref{Frepre}), we see that $T_F(\mathbf{L})=(2-F_1-F_2)(\mathbf{L})$ if and only if $\tau((\mathbf{L},\boldsymbol\infty))=0$, proving (1). The statement (2) follows from Lemma \ref{main1}, Theorem \ref{mainthm}, and the fact $F=F_1F_2$.
\end{pf} \qed

As indicated in Lemma \ref{tauform}, given a bivariate distribution function $H$ with $\mathrm{supp}(H)=[\mathbf{L},\boldsymbol\infty)\subset\mathbb{R}^2$, there exists a unique bivariate distribution function $F$ having the same support as $H$ so that
\begin{equation} \label{lambda=1}
F(\mathbf{x})=\frac{H_1(x_1)H_2(x_2)}{H_1(x_1)+H_2(x_2)-H(\mathbf{x})},\;\;\;\;\;
\mathbf{x}\in(\mathbf{L},\boldsymbol\infty).
\end{equation} Since $F$ and $H$ have the same marginals and $T_F=1-H$, by virtue of Theorem \ref{mainthm}, we have:

\begin{cor} \label{Poisson}
Let $F$ and $H$ be distribution functions on $\mathbb{R}^2$ with support $[\mathbf{L},\boldsymbol\infty)$. If $F$ and $H$ satisfy \emph{(\ref{lambda=1})}, then $F\in\mathcal{ID}(\maxconv\maxconv)$.
\end{cor}

\section{Bi-freely Max-Infinite Divisibility and Copula}
Expression (\ref{lambda=1}) suggests that bi-freely max-infinite divisibility has close ties with the following special kind of copulas
\begin{equation} \label{generalC}
C(u,v)=\frac{uv}{f(u,v)},
\end{equation} where $f$ is a function defined on $(0,1]^2$ satisfying
\begin{cond} \label{fcond}
\emph{\begin{enumerate} [$\qquad(i)$]
\item {$f(\cdot,1)\equiv1\equiv f(1,\cdot)$;}
\item {the mapping $u\mapsto f(u,v)-u$ is decreasing on $(0,1]$ for any fixed $v\in(0,1]$,
and so is the mapping $v\mapsto f(u,v)-v$ when $u$ is fixed; and}
\item {$-f$ is quasi-monotone on $(0,1]^2$.}
\end{enumerate}}
\end{cond}

Before proceeding further, let us concisely explain why the function $C$ so obtained in (\ref{generalC}) is indeed a copula. First of all, condition (\emph{ii}) implies that $C$ increases on $(0,1]^2$, while condition (\emph{iii}) with Lemma \ref{ineq} ensures that $C$ is quasi-monotone. Since $C$ meets the inequalities in (\ref{Cunifconti}) on $(0,1]^2$ by Remark \ref{qmonotone}, it admits a continuous and quasi-monotone extension to $[0,1]^2$. That the extension satisfies boundary conditions $C(\cdot,0)\equiv0\equiv C(0,\cdot)$ is due to the fact $f(u,v)\geq\max\{u,v\}$.

One apparent example of satisfying the conditions listed in Condition \ref{fcond} is the constant function $f\equiv1$, which gives the independence copula $\Pi=uv$ in (\ref{generalC}). Except this, the necessary and sufficient conditions for a non-constant function $f$ to meet Condition \ref{fcond} is that it has the representation
\begin{equation} \label{fform1}
f(u,v)=1-\theta^{-1}D(\theta\overline{u},\theta\overline{v}),\;\;\;\;\;(u,v)\in(0,1]^2,
\end{equation}
for some $\theta\in(0,1]$ and copula $D$.
Using the survival copula instead provides another way of expression:
\begin{equation} \label{fform2}
f(u,v)=\theta^{-1}-1+u+v-\theta^{-1}\widehat{D}(1-\theta\overline{u},1-\theta\overline{v}).
\end{equation}

Returning to the aforementioned conclusion, only the necessity requires a proof. First, Condition \ref{fcond} says that $f$ has an extension to $[0,1]^2$ satisfying $0\leq f(u_2,v_2)-f(u_1,v_1)\leq u_2-u_1+v_2-v_1$ for $(u_1,v_1)\leq(u_2,v_2)$. If $\theta\in(0,1)$ is any number so that $\theta(1+f(0,0))\leq1$, then the function $D:([0,\theta]\cup\{1\})^2\to[0,1]$ defined as
\[D(u,v)=\theta-\theta f(1-\theta^{-1}u,1-\theta^{-1}v),\;\;\;\;\;(u,v)\in[0,\theta]^2,\] and $D(u,1)\equiv u$ and $D(1,v)\equiv v$ is a \emph{subcopula} \cite{Copula}. Then the fact that any subcopula extends as a copula yields the representation (\ref{fform1}).

The attentive reader can notice that in addition to Condition \ref{fcond}, if $f$ satisfies the boundary conditions $f(u,0)\equiv u$ and $f(0,v)\equiv v$, then $\theta=1$ can be selected in (\ref{fform1}).

We are now in a position to state the main theorem of this section.

\begin{thm} \label{BMIDcopula1}
Let $F$ be a bivariate distribution function with marginals $F_1$ and $F_2$ having support bounded from below. Then $F$ is bi-freely max-i.d. if and only if it is of the form
\begin{equation} \label{fform3}
F=\frac{\theta F_1F_2}{\theta-D(\theta\overline{F_1},\theta\overline{F_2})},
\end{equation} where $D$ is a copula and $\theta\in(0,1]$.
\end{thm}

\begin{pf} There is nothing to prove if $F$ is the distribution function of an a.s. constant random vector, so we avoid this trivial case throughout the proof.

Suppose that $F$ is bi-freely max-i.d., and express $T_F=\lambda(1-H)$ on $\{F>0\}$ according to Theorem \ref{mainthm}. Since $\lambda>0$, we may assume $\lambda\geq1$. Write $H=C(H_1,H_2)$ for some copula $C$, and define $f(u,v)=1-\theta^{-1}D(\theta\overline{u},\theta\overline{v})$ on $[0,1]^2$, where $D=\widehat{C}$ and $\theta=1/\lambda$.
Since $\overline{H_j}=\theta\overline{F_j}$ for $j=1,2$, it follows that on the set $\{F>0\}$, we have
\[f(F_1,F_2)=1-\theta^{-1}D(\overline{H_1},\overline{H_2})=1-\theta^{-1}\overline{H}=\frac{F_1F_2}{F}.\] This confirms the ``only if'' part.

Conversely, let $F=F_1F_2/f(F_1,F_2)$, where $f$ is as in (\ref{fform1}).
Then we get $\{F>0\}=\{F_1>0\}\times\{F_2>0\}$ because $f>0$ on $(0,1]^2$, and $Q_F=1-\theta^{-1}D(\theta\overline{F_1},\theta\overline{F_2})$ holds on $\{F>0\}$. Applying (\ref{Cunifconti}) to the copula $D$, one can see that $-Q_F$ is quasi-monotone and satisfies (\ref{QFprop}) on $\overline{\{F>0\}}$. Hence the ``if'' part is derived from Theorem \ref{mainthm}.
\end{pf} \qed

Thanks to Theorem \ref{BMIDcopula1}, examining if an $F$ belongs to the set $\mathcal{ID}(\maxconv\maxconv)$ is equivalent to checking whether the underlying copula describing the dependence structure between the marginals of $F$ is of the form (\ref{generalC}). To put it differently, merely a verification of whether the bivariate function $uv/C(u,v)$ meets Condition \ref{fcond} is in demand. Below we present several copulas of this type.

Choosing the comonotone copula $D=\min\{u,v\}$ in (\ref{fform3}) gives
\begin{equation} \label{min}
\min\{u,v\}.
\end{equation} Using the survival copula in (\ref{fform3}) provides the parametric family of copulas
\[\frac{\theta uv}{1-\theta+\theta u+\theta v-D(1-\theta\overline{u},1-\theta\overline{v})},\;\;\;\;\;\theta\in(0,1].\] In particular, for any copula $D$, the copula
\begin{equation} \label{first}
\frac{uv}{u+v-D(u,v)}
\end{equation}
is of the form (\ref{generalC}), which has been proposed in (\ref{lambda=1}).

Suppose that the mixed second order partial derivative of $C$ is continuous. Then translating Condition \ref{fcond} into conditions
\begin{equation} \label{partialD}
0\leq\frac{\partial f}{\partial u},\frac{\partial f}{\partial v}\leq1\;\;\;\;\;\mathrm{and}\;\;\;\;\;
\frac{\partial^2f}{\partial u\partial v}\leq0
\end{equation} on the partial derivatives of $f(u,v)=uv/C(u,v)$ on $(0,1)^2$ yields the following result.

\begin{cor} \label{generalC2}
Let $C$ be a $C^2$-copula, and let $F_1$ and $F_2$ be distribution functions on $\mathbb{R}$ with support bounded from below. Then $C(F_1,F_2)$ is bi-freely max-i.d. if and only if the function $f=uv/C(u,v)$ satisfies \emph{(\ref{partialD})} on $(0,1)^2$.
\end{cor}

There are numerous parametric families of copulas, which have been extensively studied in the literature and also conform to the type (\ref{generalC}). For instance, one family of interest is the Lomax copulas \cite{Joenew}:
\begin{equation} \label{Lomax}
C_{p,\theta}(u,v)=\frac{uv}{\big[1-\theta(1-u^{1/p})(1-v^{1/p})\big]^p},
\end{equation} where the parameters $p>0$ and $-p\leq\theta\leq1$.
By checking the conditions in (\ref{partialD}), it is easy to see that the Lomax copula $C_{p,\theta}$ is of the form (\ref{generalC}) if and only if $\theta\in[0,1]$ and $p\in(0,1]$.

Letting $\theta=1$, the copulas in (\ref{Lomax}) become
\[(u^{-1/p}+v^{-1/p}-1)^{-p},\;\;\;\;\;p>0,\]
which is nothing but the family of Clayton copulas. Thus, a parametric Clayton copula belongs to the type (\ref{generalC}) if and only if $p\in(0,1]$.

Another well-known collection of parametric copulas is the Ali-Mikhail-Haq family:
\[\frac{uv}{1-\theta(1-u)(1-v)},\;\;\;\;\;\theta\in[-1,1],\] which can be obtained by simply choosing $p=1$ in (\ref{Lomax}). Once again, in order for copulas in this family to belong to the type (\ref{generalC}), we require $\theta\in[0,1]$.

One can employ (\ref{partialD}) to prove that the family of Farlie-Gumbel-Morgenstern copulas
\[uv\big[1+\theta(1-u)(1-v)\big],\;\;\;\;\;\theta\in[-1,1],\] has the form in (\ref{generalC}) if and only if $\theta\in[0,1]$.

\begin{remark} \emph{It is worth mentioning that if $f$ meets Condition \ref{fcond}, then so does the function $(u,v)\mapsto f^p(u^{1/p},v^{1/p})$ for any $0<p\leq1$. This is due to the inequality
\[\int_a^bt^{p-1}\;dt\geq\int_{a'}^{b'}t^{p-1}\;dt\] for any points $a,a',b,b'\in[0,\infty)$ so that $0\leq b'-b\leq a'-a$. Therefore,
\[(u,v)\mapsto\frac{uv}{f^p(u^{1/p},v^{1/p})}\] is again a copula of the type (\ref{generalC}).}
\end{remark}

We end this section with copulas in (\ref{generalC}) consisting of Pickands dependence measures. For more details regarding the copulas of this sort, we refer the reader to Section 8, where the bi-free extreme value theory is studied.

\begin{exam} \label{extremeC2}
\emph{(Extreme-value Copula). Let $\rho$ be a Pickands dependence measure on the simplex $\{(x,y)\in[0,1]^2:x+y=1\}$. It is quite straightforward to verify that the function
\[f_\rho(u,v)=1-\int\min\big\{(1-u)x,(1-v)y\big\}\;d\rho\] meets (\emph{i})-(\emph{iii}) of Condition \ref{fcond}. Letting
$C_\rho^*$ be the bivariate extreme-value copula (\ref{extremeC1}) determined by $\rho$, one can alternatively represent $f$ as
\[f_\rho(u,v)=-1+u+v-\log C_\rho^*\big(e^{-(1-u)},e^{-(1-v)}\big).\] Using the fact $\lim_{n\to\infty}n(c^{1/n}-1)=\log c$ for any constant $c>0$ and the continuity of copulas yields that for any $(u,v)\in(0,1]^2$,
\[\lim_{n\to\infty}f^n(u^{1/n},v^{1/n})=\frac{uv}{C_\rho^*(u,v)}.\] This result consequently shows that the copula $uv/f_\rho(u,v)$ lies in the max-domain of attraction of $C_\rho^*$.}
\end{exam}

\section{Example} A bivariate distribution is said to be \emph{full} if its support is neither at a point nor on a straight line, otherwise it is called \emph{non-full}. The bi-freely max-infinite divisibility of a non-full planar probability measure is characterized below.

\begin{thm} \label{nonfull}
Let $\mu$ be a Borel probability measure with support on the line $ax_1+bx_2+c=0$, but not a Dirac measure on $\mathbb{R}^2$. Then $\mu$ is bi-freely max-i.d. if and only if
\begin{enumerate} [$\qquad(1)$]
\item {the support of $\mu$ is bounded from below, and}
\item {either $b=0$ or $-a/b\in[0,\infty)$.}
\end{enumerate} In such a situation, the distribution function $F$ of $\mu$ is given by $F=\min\{F_1,F_2\}$.
\end{thm}

\begin{pf} If $a=0$ or $b=0$, then $\mu$ is the product measure of its marginals. In either situation, Corollary \ref{prod}(2) shows that $\mu$ is bi-freely max-i.d. if and only if the support of $\mu$ is bounded from below.

In the remaining proof, we will assume that $a=1$ and $b\neq0$ by multiplying the supporting line equation by a non-zero constant if necessary. We shall also assume $c=0$ by shifting the supporting line if needed. Denote by $F$ the distribution function of $\mu$.

If $F$ is bi-freely max-i.d., then by Lemma \ref{main1}, (1) is established. The relation  $\{F>0\}=\{F_1>0\}\times\{F_2>0\}$ then leads to that $\mathbf{x}\wedge\mathbf{y}\in\{F>0\}$ whenever $\mathbf{x},\mathbf{y}\in\{F>0\}$. Hence we must have $b<0$.

Conversely, suppose that $\inf\mathrm{supp}(\mu)=\mathbf{L}>-\infty$ and $b<0$. Then $\overline{\{F_j>0\}}=[L_j,\infty)$ for $j=1,2$. Pick any point
$\mathbf{x}=(x_1,x_2)$ in $\mathbb{R}^2$. If $x_1+bx_2\geq0$, then $F_1(x_1)=F(x_1,-x_1/b)$ and $F_2(x_2)=F(\mathbf{x})$ show that $F=\min\{F_1,F_2\}$ at $\mathbf{x}$. Similar reasoning handles the case $x_1+bx_2\leq0$ and also gives $F=\min\{F_1,F_2\}$ at $\mathbf{x}$. These discussions result that $F=\min\{F_1,F_2\}$, which is bi-freely max-i.d. according to (\ref{min}).
\end{pf} \qed

Denote by $\gamma_{(\mathbf{v},\mathbf{A},0)}$ the \emph{bi-free Gaussian distribution} with bi-free L\'{e}vy triplet $(\mathbf{v},\mathbf{A},0)$ \cite{HHW}. The distribution $\gamma_{(\mathbf{v},\mathbf{A},0)}$ is non-full if and only if the matrix $\mathbf{A}$ is singular, in which case $\gamma_{(\mathbf{v},\mathbf{A},0)}$ is supported on the line $\langle\mathbf{u},\mathbf{x}\rangle=\langle\mathbf{u},\mathbf{v}\rangle$, where $\mathbf{u}$ is a non-zero vector in the kernel of $\mathbf{A}$. By virtue of Theorem \ref{nonfull}, we have:

\begin{cor} \label{BFGaussin1}
The distribution $\gamma_{(\mathbf{v},\mathbf{A},0)}$ is bi-freely max-i.d. if $\mathbf{A}$ is singular and the off-diagonal element of $\mathbf{A}$ is non-negative.
\end{cor}

Next, we consider the bi-free Gaussian distribution $\gamma_c$ with correlation coefficient $c$, i.e., $\gamma_c$ has the bi-free L\'{e}vy triplet $(\mathbf{0},\mathbf{A},0)$, where
\[\mathbf{A}=\left(
\begin{array}{cc}
1 & c \\
c & 1 \\
\end{array}
\right)\] and $|c|\leq1$. We have seen that $\gamma_c$ is non-full if and only if $c=\pm1$. When $|c|<1$, $\gamma_c$ is absolutely continuous with respect to the Lebesgue measure, and its density is given by the formula \cite{HW1}:
\[p_c(s,t)=\frac{1-c^2}{4\pi^2}\frac{\sqrt{4-s^2}\sqrt{4-t^2}}{D_c(s,t)},\;\;\;\;\;(s,t)\in[-2,2]^2,\] where
\[D_c(s,t)=(1-c^2)^2-c(1+c^2)st+c^2(s^2+t^2).\] Notice that $D_c(s,t)=(2cs-(1+c^2)t)^2/4+(1-c^2)^2(4-t^2)/4>0$ on $[-2,2]^2$.

Before we prove that $\gamma_c$ is not bi-freely max-i.d. except for $c=0,1$, two lemmas are required for the verification.

\begin{lem} \label{identity}
Let $-1<c<1$. For any $x\in[-2,2]$, we have
\[\int_{-2}^2\frac{\sqrt{4-t^2}}{D_c(x,t)}\;dt=\frac{2\pi}{1-c^2}.\]
\end{lem}

\begin{pf} Let $(X,Y)$ be a random vector having distribution $\gamma_c$ with $c\neq\pm1$. Since $X$ has probability density function $\sqrt{4-x^2}/(2\pi)$, the semi-circular law, it follows that
\[\int_{-2}^x\left(\int_{-2}^2p_c(s,t)\;dt\right)ds=\mathrm{Pr}[-2\leq X\leq x]=\frac{1}{2\pi}\int_{-2}^x\sqrt{4-s^2}\;ds\] for any $x\in[-2,2]$. Differentiating this equation with respect to $x$ gives that
\[\int_{-2}^2p_c(x,t)\;dt=\frac{\sqrt{4-x^2}}{2\pi},\;\;\;\;\;x\in[-2,2],\] inferring the desired result after some algebraic works.
\end{pf} \qed

\begin{lem} \label{compare}
Let $0<c<1$. Then we have
\begin{equation} \label{compare1}
\int_{-2}^x\sqrt{4-s^2}\left[\int_{-2}^y\left(\frac{\sqrt{4-t^2}}{D_{-c}(s,t)}-\frac{\sqrt{4-t^2}}{D_{-c}(x,t)}\right)dt\right]ds\;<0
\end{equation}
for any $(x,y)\in(-2,2)^2$ and
\begin{equation} \label{compare2}
\lim_{x\downarrow-2}\int_{-2}^x
\frac{\sqrt{4-s^2}}{D_c(x,-2)}ds\bigg/
\int_{-2}^x\frac{\sqrt{4-s^2}}{D_c(s,-2)}\;ds=1.
\end{equation}
\end{lem}

\begin{pf} The desired inequality (\ref{compare1}) follows from the inequality
\begin{equation} \label{compare3}
\int_{-2}^x\sqrt{4-s^2}\left[\int_{-2}^y\left(\frac{\sqrt{4-t^2}}{D_c(s,t)}-\frac{\sqrt{4-t^2}}
{D_c(x,t)}\right)dt\right]ds\;>0,\;\;\;(x,y)\in(-2,2)^2.
\end{equation} Indeed, with the help of Lemma \ref{identity} and the relation $D_{-c}(s,t)=D_c(-s,t)$, one can see that for any $(x,y)\in(-2,2)^2$,
\begin{align*}
\int_{-2}^x\int_{-2}^y\frac{\sqrt{4-s^2}\sqrt{4-t^2}}{D_{-c}(s,t)}\;dsdt&=
\int_{-x}^2\int_{-2}^y\frac{\sqrt{4-s^2}\sqrt{4-t^2}}{D_c(s,t)}\;dsdt \\
&=
\frac{4\pi^2F_2(y)}{1-c^2}-\int_{-2}^{-x}\int_{-2}^y\frac{\sqrt{4-s^2}\sqrt{4-t^2}}{D_c(s,t)}\;dsdt \\
&<\frac{4\pi^2F_2(y)}{1-c^2}-\int_{-2}^{-x}\int_{-2}^y\frac{\sqrt{4-s^2}\sqrt{4-t^2}}{D_c(-x,t)}\;dsdt \\
&=\int_{-x}^2\int_{-2}^y\frac{\sqrt{4-s^2}\sqrt{4-t^2}}{D_c(-x,t)}\;dsdt \\
&=\int_{-2}^x\int_{-2}^y\frac{\sqrt{4-s^2}\sqrt{4-t^2}}{D_{-c}(x,t)}\;dsdt.
\end{align*}

From now to the end of the proof, $x\in(-2,2)$ is a fixed point. Then
\[y_0=\frac{2cx}{1+c^2}\in(-2,2).\] The inference strategy of proving the inequality (\ref{compare3}) is to distinguish two possibilities $y\in[y_0,2)$ and $y\in(-2,y_0)$.

First, consider the function
\[f(s,y)=\int_{-2}^y\left(\frac{\sqrt{4-t^2}}{D_c(s,t)}-\frac{\sqrt{4-t^2}}{D_c(x,t)}\right)dt\] defined on $[-2,x)\times[y_0,2]$. According to Lemma \ref{identity}, we have $f(s,2)=0$. Since
\[\frac{\partial}{\partial y}f(s,y)=c(x-s)\sqrt{4-y^2}\cdot\frac{c(x+s)-(1+c^2)y}{D_c(x,y)D_c(s,y)}\] and
$c(x+s)-(1+c^2)y<2cx-(1+c^2)y_0=0$ on $[-2,x)\times[y_0,2)$, it follows that
$f(s,y)>f(s,2)=0$ for any point $(s,y)\in[-2,x)\times[y_0,2)$. Consequently, (\ref{compare3}) is true for $(x,y)$ with $y\in[y_0,2)$.

Next, let $y\in(-2,y_0)$ be an arbitrary but fixed point. Then
\[\xi_0=\max\left\{\frac{(1+c^2)y}{2c},-2\right\}\in[-2,2).\] Further, let
\[g(\xi,t)=\int_{-2}^\xi\left(\frac{\sqrt{4-s^2}}{D_c(s,t)}-\frac{\sqrt{4-s^2}}{D_c(\xi,t)}\right)ds,\;\;\;\;\;
(\xi,t)\in[\xi_0,2]\times[-2,y).\] Since
\[\frac{\partial}{\partial\xi}g(\xi,t)=\frac{c\big[2c\xi-(1+c^2)t\big]}{D_c^2(\xi,t)}\int_{-2}^\xi\sqrt{4-s^2}\;ds\] and
$2c\xi-(1+c^2)t\geq2c\xi_0-(1+c^2)t>0$ on $[\xi_0,2]\times[-2,y)$, it follows that $g(\xi,t)>g(\xi_0,t)$ for any $(\xi,t)\in(\xi_0,2)\times[-2,y)$. Using the fact $x>\xi_0$, we derive that
\begin{align*}
\int_{-2}^y\sqrt{4-t^2}g(x,t)\;dt&>\int_{-2}^y\sqrt{4-t^2}g(\xi_0,t)\;dt \\
&=\int_{-2}^{\xi_0}\sqrt{4-s^2}\left[\int_{-2}^y\left(\frac{\sqrt{4-t^2}}{D_c(s,t)}
-\frac{\sqrt{4-t^2}}{D_c(\xi_0,t)}\right)dt\right]ds.
\end{align*} The last double integral displayed above is obviously zero if $\xi_0=-2$, while it is strictly positive if $\xi_0>-2$ by virtue of the established result in the first case. Hence (\ref{compare3}) is also valid for $(x,y)$ with $y\in(-2,y_0)$.

Proving (\ref{compare2}) is equivalent to computing the limit
\[\lim_{x\downarrow-2}\int_{-2}^x\left[\frac{\sqrt{4-s^2}}{D_c(x,-2)}-
\frac{\sqrt{4-s^2}}{D_c(s,-2)}\right]ds\bigg/
\int_{-2}^x\frac{\sqrt{4-s^2}}{D_c(s,-2)}\;ds,\] which can be done using the L'H\^{o}pital's rule:
\begin{align*}
-\frac{\partial D_c/\partial x (-2,-2)}{D_c(-2,-2)}\lim_{x\downarrow-2}\int_{-2}^x\sqrt{4-s^2}ds\bigg/\sqrt{4-x^2}=\frac{2c}{(1-c)^2}
\lim_{x\downarrow-2}\frac{4-x^2}{x}=0.
\end{align*} Hence the proof is complete.
\end{pf} \qed

\begin{thm} \label{example}
The bi-free Gaussian distribution $\gamma_c$ of correlation coefficient $c$ is bi-freely max-i.d. if and only if either $c=0$ or $c=1$. Furthermore, the copulas associated with $\gamma_0$ and $\gamma_1$ are respectively the dependence copula and comonotone copula.
\end{thm}

\begin{pf} That $\gamma_c$ is bi-freely max-i.d. for the cases $c=0$ and $c=1$ have been respectively verified in Corollary \ref{prod}(2) and Corollary \ref{BFGaussin1}. The situation $c=-1$ also has been discussed in Corollary \ref{BFGaussin1}.

Now, assume $c\in(-1,0)$ and denote by $F$ the distribution function of $\gamma_c$. We shall prove that the mapping $x\mapsto Q_F(x,y)$ strictly decreases on $(-2,2)$ for any fixed $y\in(-2,2)$, which together with Remark \ref{QFremark} yields that $F$ is not bi-freely max-infinitely divisible. Passing to partial derivatives, this is equivalent to showing that $F_1'(x)F(x,y)-F_1(x)F_x(x,y)<0$ on $(-2,2)^2$. One can derive this desired inequality by making use of Lemma \ref{compare} and the following two identities
\[F_1'(x)F(x,y)=\frac{1-c^2}{4\pi^2}\frac{\sqrt{4-x^2}}{2\pi}\int_{-2}^x
\sqrt{4-s^2}\left(\int_{-2}^y\frac{\sqrt{4-t^2}}{D_c(s,t)}\;dt\right)ds\]
and
\[F_1(x)F_x(x,y)=\frac{1-c^2}{4\pi^2}\frac{\sqrt{4-x^2}}{2\pi}\int_{-2}^x
\sqrt{4-s^2}\left(\int_{-2}^y\frac{\sqrt{4-t^2}}{D_c(x,t)}\;dt\right)ds.\]

Finally, only the situation of $c\in(0,1)$ remains to process. Our inference strategy for this case is to claim that $\partial T_F/\partial x>0$ at some point in $(-2,2)^2$ near point $(-2,-2)$. Once this claim is accomplished, we have shown that the mapping $x\mapsto T_F(x,y_0)$ strictly increases over some open interval containing $x_0$. Therefore, $\gamma_c$ cannot be bi-freely max-i.d. according to Lemma \ref{main1}.

For this, observe first that we have
\begin{align*}
\frac{\partial}{\partial x}T_F(x,y)&=\frac{F_2(y)}{F(x,y)}\cdot\frac{F_1'(x)F(x,y)-F_1(x)F_x(x,y)}{F(x,y)}-F_1'(x) \\
&=\frac{\sqrt{4-x^2}}{2\pi}
\left[\frac{I(x,y)F_2(y)}{F(x,y)}-1\right],
\end{align*} where the function $I(x,y)$ is given by
\[\int_{-2}^x\int_{-2}^y\sqrt{4-s^2}\left[\frac{\sqrt{4-t^2}}{D_c(s,t)}-
\frac{\sqrt{4-t^2}}{D_c(x,t)}\right]dtds\bigg/\int_{-2}^x\int_{-2}^y\frac{\sqrt{4-s^2}\sqrt{4-t^2}}
{D_c(s,t)}\;dsdt.\] What will be established is that
\begin{equation} \label{infty}
\lim_{x\downarrow-2}\lim_{y\downarrow-2}\frac{I(x,y)F_2(y)}{F(x,y)}=\infty,
\end{equation} which completes the proof of the claim.

The limit (\ref{infty}) is merely an application of the L'H\^{o}pital's rule. Indeed, we have
\begin{align*}
&\lim_{y\downarrow-2}I(x,y) \\
=&\lim_{y\downarrow-2}\int_{-2}^x\sqrt{4-s^2}\left[\frac{\sqrt{4-y^2}}{D_c(s,y)}-
\frac{\sqrt{4-y^2}}{D_c(x,y)}\right]ds\bigg/\int_{-2}^x\frac{\sqrt{4-s^2}\sqrt{4-y^2}}{D_c(s,y)}\;ds \\
=&\int_{-2}^x\left[\frac{\sqrt{4-s^2}}{D_c(s,-2)}-
\frac{\sqrt{4-s^2}}{D_c(x,-2)}\right]ds\bigg/\int_{-2}^x\frac{\sqrt{4-s^2}}{D_c(s,-2)}\;ds
\end{align*} and
\[\lim_{y\downarrow-2}\frac{F_2(y)}{F(x,y)}=2\pi\bigg/\left((1-c^2)\int_{-2}^x\frac{\sqrt{4-s^2}}{D_c(s,-2)}\;ds\right).\] These imply that the limit (\ref{infty}) simplifies into
\[\frac{2\pi}{1-c^2}\lim_{x\downarrow-2}\int_{-2}^x\left[\frac{\sqrt{4-s^2}}{D_c(s,-2)}-
\frac{\sqrt{4-s^2}}{D_c(x,-2)}\right]ds\bigg/\left(\int_{-2}^x\frac{\sqrt{4-s^2}}{D_c(s,-2)}\;ds\right)^2.\] Thanks to (\ref{compare2}), the limit displayed above is further equal to
\begin{align*}
&\lim_{x\downarrow-2}D_c'(x,-2)\cdot
\int_{-2}^x\frac{\sqrt{4-s^2}}{D_c(x,-2)}\;ds\bigg/\left(2\sqrt{4-x^2}\int_{-2}^x\frac{\sqrt{4-s^2}}{D_c(s,-2)}\;ds\right) \\
=&\lim_{x\downarrow-2}D_c'(x,-2)\bigg/\left(2\sqrt{4-x^2}\right) \\
=&\lim_{x\downarrow-2}\frac{c(1-c)^2+c^2(x+2)}{\sqrt{4-x^2}}=\infty,
\end{align*} which is exactly what we want to gain. Hence $\gamma_c$ is not bi-freely max-i.d. if $c\in(0,1)$.
\end{pf} \qed

The last example of this section is the \emph{bi-free compound Poisson distribution} $P_{(\lambda,\tau)}$ with rate $\lambda>0$ and jump distribution $\tau$. It was shown in \cite{HHW} that $P_{(\lambda,\tau)}$ is non-full if and only if $\tau$ is non-full, in which case both of $P_{(\lambda,\tau)}$ and $\tau$ are supported on the same line in the plane. Using Theorem \ref{nonfull} again infers the following corollary.

\begin{cor} \label{BFPoisson}
The distribution $P_{(\lambda,\tau)}$ is bi-freely max-i.d. if the support of $\tau$ lies on a line with a positive slope or on a vertical line. In particular, $P_{(\lambda,\delta_{(\alpha,\beta)})}$ is bi-freely max-i.d. if and only if $\alpha\beta\geq0$.
\end{cor}

\section{Max-Infinite Divisibility} We now turn our attention to exploring the connection between the set $\mathcal{ID}(\maxconv\maxconv)$ and the family of max-i.d. distribution functions on $\mathbb{R}^2$. Recall that a max-i.d. distribution function is characterized by (\ref{exponent}) and it may not have support bounded from below. Our first result is that $\mathcal{ID}(\maxconv\maxconv)$ is a strict subset of max-i.d. distribution functions.

\begin{thm} \label{classical1}
Any bi-freely max-infinitely divisible distribution function is max-infinitely divisible.
\end{thm}

\begin{pf} Let $F$ be bi-freely max-infinitely divisible. Because the distribution function of an a.s. constant random vector is automatically max-i.d., we may assume that the $F$ being considered is not such a case. To attain the max-infinite divisibility of $F$, it suffices to show that $F^{1/n}$ is quasi-monotone on $\mathbb{R}^2$ for any $n\in\mathbb{N}$.

As recorded in Lemma \ref{tauform} and Theorem \ref{mainthm}, $\log F$ is quasi-monotone on $\{F>0\}$. Then it follows from Lemma \ref{ineq} that for any points $\mathbf{x}\leq\mathbf{y}$ in $\{F>0\}$,
\[V_{F^{1/n}}([\mathbf{x},\mathbf{y}])=e^{n^{-1}\log F(\mathbf{y})}-e^{n^{-1}\log F(x_1,y_2)}-
e^{n^{-1}\log F(y_1,x_2)}+e^{n^{-1}\log F(\mathbf{x})}\geq0.\] This confirms that $F^{1/n}$ is quasi-monotone on $\{F>0\}$, and consequently on $\mathbb{R}^2$ as well. This finishes the proof.
\end{pf} \qed

Due to the theorem above, we can apply the known conclusions regarding max-infinite divisibility in the literature to the bi-free case, such as the following one from \cite{Maxclassical}.

\begin{cor} \label{support1}
The support $S$ of any $\mu\in\mathcal{ID}(\maxconv\maxconv)$ is characterized by the properties that $\mathbf{x}\vee\mathbf{y}\in S$ whenever $\mathbf{x},\mathbf{y}\in S$ and $S$ contains a sequence $\mathbf{x}_n\to\inf S$.
\end{cor}

There are more can be said about the max-infinite divisibility concept in classical and bi-free probability.

\begin{thm} \label{classical2}
Given an $F\in\mathcal{ID}(\maxconv\maxconv)$,
the function
\begin{equation} \label{maxGTF}
G(\mathbf{x})=\left\{
\begin{array}{ll}
\exp[-t T_F(\mathbf{x})], & \hbox{$\mathbf{x}\in[\mathbf{L},\boldsymbol\infty)$}, \\
0, & \hbox{otherwise,}
\end{array}\right.
\end{equation} is a max-i.d. distribution function for any $t>0$, where $\mathbf{L}=\inf\mathrm{supp}(F)$.

Conversely, a max-i.d. distribution function $G$ satisfying $\{G>0\}=[\mathbf{L},\boldsymbol\infty)$ for some $\mathbf{L}\in\mathbb{R}^2$ is of the form \emph{(\ref{maxGTF})} for some $t>0$ and some $F\in\mathcal{ID}(\maxconv\maxconv)$ with $\inf\mathrm{supp}(F)=\mathbf{L}$.
\end{thm}

\begin{pf} Let $F\in\mathcal{ID}(\maxconv\maxconv)$. It is clear that $G$ is max-i.d. if
$F$ is the distribution function of the Dirac measure $\delta_\mathbf{L}$, because $T_F=0$ on $[\mathbf{L},\boldsymbol\infty)$ in this case. Now, suppose that $F$ is not the previously mentioned case, and let $t>0$. Then by Lemma \ref{main1}, the function $G$ defined in (\ref{maxGTF}) is a bivariate distribution function. Indeed, Lemma \ref{ineq} and the quasi-monotonicity of $-tT_F$ yield at once that $G$ is also quasi-monotone on $\mathbb{R}^2$. Then defining $G^{(n)}=\exp[-tT_F/n]$ on $[\mathbf{L},\boldsymbol\infty)$ and zero elsewhere results in $(G^{(n)})^n=G$, proving that $G$ is max-infinitely divisible.

Conversely, suppose that $G$ is max-i.d. and $G>0$ on $E=[\mathbf{L},\boldsymbol\infty)\subset\mathbb{R}^2$. Let $\widetilde{\tau}$ be an exponent measure of $G$, i.e., one represents $G$ as
\[G(\mathbf{x})=\exp\big\{-\widetilde{\tau}\big(E\backslash(-\boldsymbol\infty,\mathbf{x}]\big)\big\},
\;\;\;\;\;\mathbf{x}\in E,\] by (\ref{exponent}). Then according to the hypothesis, we have $\widetilde{\tau}(E\backslash\{\mathbf{L}\})=-\log G(\mathbf{L})<\infty$. Next, pick a finite number $t>\widetilde{\tau}(E\backslash\{\mathbf{L}\})$ and define $\tau=t^{-1}\mathbf{1}_{E\backslash\{\mathbf{L}\}}\widetilde{\tau}$, where $\mathbf{1}_{E\backslash\{\mathbf{L}\}}$ is the indicator function of $E\backslash\{\mathbf{L}\}$. With this $\tau$, one can construct distribution functions $F_j$ and $F$ with the stated properties in Lemma \ref{tauform}. Such an $F$ is bi-freely max-i.d. by Theorem \ref{mainthm}, and it is fairly easy to verify that $G$ is written as (\ref{maxGTF}). This finishes the proof.
\end{pf} \qed

\section{Bi-free Extreme value Theory}
Recall that a bivariate extreme value distribution can be conveniently expressed by using its marginals, both of which are univariate extreme value distribution functions, and its Pickands dependence function. In this section, we shall provide the bi-free extreme value results and discuss their relevance to classical extreme value theory. We first begin with one definition from \cite{V15EXT}.

\begin{pdef} \emph{Let $F$ be a bivariate distribution function. If for any $n\in\mathbb{N}$, there are $a_n>0$ and $c_n>0$ and real numbers $b_n$ and $d_n$ so that
\[F^{\maxconvexp\maxconvexp n}(a_nx_1+b_n,c_nx_2+d_n)=F(x_1,x_2)\] on $\mathbb{R}^2$, then $F$ is termed \emph{bi-free max-stable}.
If there are a bivariate distribution function $H$ and real numbers $a_n,b_n,c_n$ and $d_n$ with $a_n>0$ and $c_n>0$ so that
\begin{equation} \label{bif-freeD}
H^{\maxconvexp\maxconvexp n}(a_nx_1+b_n,c_nx_2+d_n)\to F(x_1,x_2)
\end{equation} weakly,
then we say that $F$ is a \emph{bi-free extreme value distribution function} and $H$
is in the \emph{bi-free max-domain of attraction} of $F$. The collection of those $H$ satisfying (\ref{bif-freeD}) is denoted by $\mathcal{D}_{\maxconvlow\maxconvlow}(F)$.}
\end{pdef}

It follows at once from the definition that both marginals of a bi-free max-stable distribution function are freely max-stable \cite{freeextre}.

As demonstrated in Example \ref{extremeC2}, a Pickands dependence function $A$ generates the copula
\begin{equation} \label{extremeC3}
C_A(u,v)=\frac{uv}{-1+u+v+(2-u-v)A\left(\frac{1-u}{2-u-v}\right)},\;\;\;\;\;(u,v)\in(0,1)^2,
\end{equation} and the copula thus obtained belongs to the domain of attraction of the extreme-value copula $C_A^*$ determined by $A$. It is shown below that same as the classical case, a bi-free extreme value distribution function is likewise determined by its marginals and a copula involving a unique dependence function.

\begin{thm} \label{BFstable}
Let $F$ be a bivariate distribution function with both marginals freely max-stable. The followings are equivalent:
\begin{enumerate} [$\qquad(1)$]
\item {$F$ is bi-free max-stable;}
\item {$F\in\mathcal{D}_{\maxconvlow\maxconvlow}(F)$;}
\item {$F$ is a bi-free extreme value distribution;}
\item {there exists some Pickands dependence function $A$ so that
\begin{equation} \label{stableform}
F=C_A(F_1,F_2),
\end{equation} where $C_A$ is the copula in \emph{(\ref{extremeC3})}.}
\end{enumerate} Moreover, if one of \emph{(1)-(4)} holds, then the dependence function $A$ in \emph{(4)} is unique.
\end{thm}

\begin{pf} The implications $(1)\Rightarrow(2)\Rightarrow(3)$ are clear.

$(3)\Rightarrow(4)$: Suppose that $H\in\mathcal{D}_{\maxconvlow\maxconvlow}(F)$, i.e., for any $n\in\mathbb{N}$, there are $a_n,c_n>0$ and $b_n,d_n\in\mathbb{R}$ so that
$H^{\maxconvexp\maxconvexp n}(a_nx_1+b_n,c_nx_2+d_n)\to F(x_1,x_2)$ weakly. For notational convenience, in the following arguments we shall write $x_{1n}=a_nx_1+b_n$ and $x_{2n}=c_nx_2+d_n$ whenever points $x_1$ and $x_2$ are given.

Since
\begin{equation} \label{freedomain}
H_1^{\maxconvexp n}(a_nx+b_n)\to F_1(x)\;\;\;\;\;\mathrm{and}\;\;\;\;\;
H_2^{\maxconvexp n}(c_nx+d_n)\to F_2(x)
\end{equation}
weakly, it follows that
\begin{equation} \label{cdomain}
H_1^n(a_nx+b_n)\to G_1(x)\;\;\;\;\;\mathrm{and}\;\;\;\;\;H_2^n(c_nx+d_n)\to G_2(x)
\end{equation}
weakly, where each $G_j$ is a univariate extreme value distribution function satisfying the relation $(1+\log G_j)_+=F_j$ \cite{freeextre}. This implies that for $j=1,2$, we have
\begin{equation} \label{limH_jG_j}
\lim_{n\to\infty}n\big[1-H_j(x_{jn})\big]=-\log G_j(x_j),\;\;\;\;\;x_j\in\{G_j>0\}.
\end{equation} The weak convergence of $H^{\maxconvexp\maxconvexp n}(x_{1n},x_{2n})\to F(x_1,x_2)$ also allows us to obtain the limit
\begin{equation} \label{limHF1}
n\left[\frac{H_1(x_{1n})H_2(x_{2n})}{H(x_{1n},x_{2n})}-1\right]\to\frac{F_1(x_1)F_2(x_2)}{F(x_1,x_2)}
-1
\end{equation}
for any point $(x_1,x_2)$ in the positive set $\mathcal{P}=\{F_1>0\}\times\{F_2>0\}$.
So far, we conclude from Theorem \ref{maxconv1} that $F$ is bi-freely max-infinitely divisible.

Next, observe that because of the limit $\lim_{n\to\infty}H(x_{1n},x_{2n})=1$ obtained by (\ref{limH_jG_j}), we are able to rewrite (\ref{limHF1}) as that for $(x_1,x_2)\in\mathcal{P}$,
\begin{equation} \label{limHF2}
n\big[H_1(x_{1n})H_2(x_{2n})-H(x_{1n},x_{2n})\big]\to\frac{F_1(x_1)F_2(x_2)}{F(x_1,x_2)}-1.
\end{equation} Decomposing the left-hand side of the limit in (\ref{limHF2}) into the sum
\[n\big[H_1(x_{1n})-1\big]+n\big[H_2(x_{2n})-1\big]+n\big[1-H(x_{1n},x_{2n})\big]+o(1/n)\] and using (\ref{limH_jG_j}) allow us to further attain that for $(x_1,x_2)\in\mathcal{P}$,
\begin{equation} \label{limHF3}
n\big[1-H(x_{1n},x_{2n})\big]\to\frac{F_1(x_1)F_2(x_2)}{F(x_1,x_2)}
-1-\log G_1(x_1)-\log G_2(x_2).
\end{equation}
On the other hand, due to the bi-freely max-infinite divisibility of $F$, there exists certain function $f$ satisfying (\emph{i})-(\emph{iii}) in Condition \ref{fcond} so that $F_1F_2/F=f(F_1,F_2)$ on $\mathcal{P}$. Then using the dependence relation $H=D(H_1,H_2)$ for some copula $D$ and letting $g(u,v)=1-u-v+f(u,v)$, the limit in (\ref{limHF3}) now becomes
\begin{equation} \label{limHF4}
n\big[1-D(H_1(x_{1n}),H_2(x_{2n}))\big]\to g(F_1(x_1),F_2(x_2)),\;\;\;\;\;(x_1,x_2)\in\mathcal{P}.
\end{equation}

For further discussion, we consider the subcopula
\begin{equation} \label{Cg}
C(u,v)=\exp\big[-g(1+\log u,1+\log v)\big],\;\;\;\;\;(u,v)\in[1/e,1]^2.
\end{equation}
With the help of this subcopula, (\ref{limHF4}) is equivalent to
\begin{equation} \label{limHF5}
n\big[1-D(H_1(x_{1n}),H_2(x_{2n}))\big]\to-\log C(G_1(x_1),G_2(x_2)),
\end{equation} where the points $x_j$ are selected from the set $\{G_j>1/e\}=\{F_j>0\}$.
Since any bivariate subcopula also meets the Lipschitz continuity condition (\ref{Cunifconti}) in its domain, it follows from (\ref{limH_jG_j}) and (\ref{limHF5}) that
\begin{equation} \label{limHF6}
n\left[1-D\left(1+\frac{\log G_1(x_1)}{n},1+\frac{\log G_2(x_2)}{n}\right)\right]\to-\log C(G_1(x_1),G_2(x_2))
\end{equation} on $\mathcal{P}$. Now, letting $x=-\log G_1(x_1)$ and $y=-\log G_2(x_2)$ consequently results in that
\begin{equation} \label{limHF7}
n\big[1-D(1-n^{-1}x,1-n^{-1}y)\big]\to-\log C(e^{-x},e^{-y})
\end{equation} is true for any $(x,y)\in[0,1)^2$.

After these discussions, consideration is given to handling points $x_j$ chosen from the set $\{G_j>0\}$. In this general situation, by choosing any $k\in\mathbb{N}$ with $-\log G_j(x_j)<k$ for $j=1,2$, we infer from the limits in (\ref{limH_jG_j}) and (\ref{limHF7}) and from the Lipschitz continuity of $D$ that
\begin{align*}
nk[1-H(x_{1(nk)},&x_{2(nk)})]=nk[1-D(H_1(x_{1(nk)}),H_2(x_{2(nk)}))] \\
&=nk\left[1-D\left(1+\frac{\log G_1(x_1)}{nk},1+\frac{\log G_2(x_2)}{nk}\right)\right]+o(1/n) \\
&\to-\log C^k\big(G_1^{1/k}(x_1),G_2^{1/k}(x_2)\big)\;\;\mathrm{as}\;\;n\to\infty.
\end{align*} Notice that the limit displayed above is independent of the values of $k$ as long as $-\log G_j(x_j)<k$ for $j=1,2$ because we have
\begin{align*}
-\log C^{k_1}\big(G_1^{1/k_1}(x_1),G_2^{1/k_1}(x_2)\big)&=
\lim_{n\to\infty}nk_1k_2[1-H(x_{1(nk_1k_2)},x_{2(nk_1k_2)})] \\
&=-\log C^{k_2}\big(G_1^{1/k_2}(x_1),G_2^{1/k_2}(x_2)\big).
\end{align*} This observation makes it possible to extend the domain of the function $C$ to $[0,1]^2$. More precisely, for any $(u,v)\in(0,1]^2$ and any $k\in\mathbb{N}$ so that $u,v>e^{-k}$, one can define
\[C(u,v)=C^k(u^{1/k},v^{1/k})\] and see that such an extension is unambiguous. Due to this extension, (\ref{limHF7}) is finally generalized to
\begin{equation} \label{limHF8}
n\big[1-D(1-n^{-1}x,1-n^{-1}y)\big]\to-\log C(e^{-x},e^{-y}),\;\;\;\;\;(x,y)\in[0,\infty)^2.
\end{equation} One can easily verify the quasi-monotonicity of $C$ on $[0,1]^2$ by (\ref{limHF8}) and Lemma \ref{ineq}. Consequently, in view of (\ref{extreCdef2}), we have confirmed that $C$ is an extreme-value copula.

Finally, making use of the Pickands dependence function $A$ of the extreme-value copula $C$ and the relation (\ref{Cg}) yields the expression
\begin{equation} \label{fArelation}
f(u,v)=-1+u+v+(2-u-v)A\left(\frac{1-u}{2-u-v}\right)
\end{equation} for $(u,v)\in[0,1]^2\backslash\{(1,1)\}$. This leads to the representation formula of $F$ in statement (4), and the proof of (3)$\Rightarrow$(4) is now accomplished. As for the uniqueness of the dependence function $A$, this follows from the fact that a free extreme value distribution is continuous on the real line.

$(4)\Rightarrow(1)$: Suppose that $F$ is of the form (\ref{stableform}) for some dependence function $A$. The boundary conditions on $A$ imply that we have $\{F>0\}=\{F_1>0\}\times\{F_2>0\}$, and $F=1$ when $F_1=1$ or $F_2=1$. Let $G$ be given by (\ref{Gextreme}), where each $G_j$ is the univariate extreme value distribution satisfying the relation $(1+\log G_j)_+=F_j$. Since $G$ is max-stable, it follows that for any $n\in\mathbb{N}$, there exist $a_n,c_n>0$ and $b_n,d_n\in\mathbb{R}$ so that
$G^n(a_nx_1+b_n,c_nx_2+d_n)=G(x_1,x_2)$ for any $x_1,x_2\in\mathbb{R}$.

In the rest of the proof, points $(x_1,x_2)$ are all selected from $\{0<F<1\}$ and, as before, we adopt the simplified notations $x_{1n}=a_nx_1+b_n$ and $x_{2n}=c_nx_2+d_n$. Observe that we have
\begin{align*}
F_1(x_1)=1+\log G_1(x_1)&=1+\log G_1^n(x_{1n}) \\
&=nF_1(x_{1n})-(n-1)=F_1^{\maxconvexp n}(x_{1n})
\end{align*} and that through the same reasoning $F_2^{\maxconvexp n}(c_nx_2+d_n)=F_2(x_2)$ is established. Moreover, the expression (\ref{Gextreme}) shows that
\[-\big(\overline{F_1}(x_1)+\overline{F_2}(x_2)\big)A\left(\frac{\overline{F_1}(x_1)}{\overline{F_1}(x_2)+\overline{F_2}(x_2)}\right)=
\log G(x_1,x_2)=n\log G(x_{1n},x_{2n}).\] All these findings together yield that
\begin{align*}
\frac{F_1(x_1)F_2(x_2)}{F(x_1,x_2)}-1&=\log G_1(x_1)+\log G_2(x_2)-\log G(x_1,x_2) \\
&=n\big[\log G_1(x_{1n})+\log G_2(x_{2n})-\log G(x_{1n},x_{2n})\big] \\
&=n\left[\frac{F_1(x_{1n})F_2(x_{2n})}{F(x_{1n},x_{2n})}-1\right].
\end{align*} Notice that the above result is also true when $F(x_1,x_2)=1$.
We consequently derive $F^{\maxconvexp\maxconvexp n}(x_{1n},x_{2n})=F(x_1,x_2)$ on $\{F>0\}$ in view of Remark \ref{maxremark}, proving that $F$ is bi-free max-stable. This accomplishes the proof of (4)$\Rightarrow$(1).
\end{pf} \qed

The attentive reader might have noticed from the preceding proof the following fact concerned with the bi-free max-domain of attraction.

\begin{thm} \label{domain}
Let $G_j$ be univariate extreme value distributions and let $(1+\log G_j)_+=F_j$ for $j=1,2$. Given a Pickands dependence function $A$, define $G$ and $F$ as in \emph{(\ref{Gextreme})} and \emph{(\ref{stableform})}, respectively. Then
\[\mathcal{D}_*(G)=\mathcal{D}_{\maxconvlow\maxconvlow}(F)\] and the normalizing vectors are also the same.
\end{thm}

\begin{pf} Let $H\in\mathcal{D}_{\maxconvlow\maxconvlow}(F)$. It has been shown in the proof of the implication (3)$\Rightarrow$(4) of Theorem \ref{BFstable} that there exist normalizing constants $a_n,c_n>0$ and $b_n,d_n\in\mathbb{R}$ so that
\begin{equation} \label{HdomainG}
\lim_{n\to\infty}n[1-H(a_nx_1+b_n,c_nx_2+d_n)]=-\log G(x_1,x_2)
\end{equation}
whenever $G(x_1,x_2)>0$. This implies that $H\in\mathcal{D}_*(G)$ by (\ref{EV2}).

Conversely, let $H\in\mathcal{D}_*(G)$. Thanks to (\ref{EV2}), it is easy to see that (\ref{limHF3}) is valid on $\{F>0\}$, so (\ref{limHF1}) is true as well. Since the limits in (\ref{freedomain}) are equivalent those in (\ref{cdomain}) \cite{freeextre}, we have shown that $H\in\mathcal{D}_{\maxconvlow\maxconvlow}(F)$. This finishes the proof.
\end{pf} \qed

One can make use of Theorem \ref{domain} and max-domain of attraction criteria to determine the bi-free max-domain attraction of a given bi-free max-stable distribution function. For the details about the criteria, we refer the reader to \cite{Resnick}.

There are a number of parametric families of bi-free extreme value distribution functions. Some of these which might be of interest are:

\begin{exam} \emph{(Gumbel Mixed Model). Consider the parametric copulas
\begin{equation} \label{Gumbel1}
C_\theta(u,v)=\frac{uv}{f_\theta(u,v)},
\end{equation} where $\theta\in[0,1]$ and
\[f_\theta(u,v)=1-\theta\cdot\frac{(1-u)(1-v)}{2-u-v},\;\;\;\;\;(u,v)\in[0,1)^2.\] Note that $C_\theta$ satisfies the formula (\ref{extremeC3}) with the Pickands dependence function $A_\theta(t)=\theta t^2-\theta t+1$. Hence by Theorem \ref{BFstable}, $C_\theta(F_1,F_2)$ is a bi-free extreme value distribution function provided that each $F_j$ is of free extreme type. Also, $C_\theta$ belongs to the domain of attraction of the extreme-value copula \[C_\theta^*(u,v)=uv\exp\left[-\theta\cdot\frac{\log u\log v}{\log(uv)}\right].\] The family $\{C_\theta^*\}$ is one of two general forms of parametric copulas proposed by Gumbel in the study of bivariate extremes.}
\end{exam}

\begin{exam} \emph{(Logistic Model). The dependence function $A_m(t)=(t^m+(1-t)^m)^{1/m}$, $m\geq1$, generates the copula
\begin{equation} \label{Gumbel2}
C_m(u,v)=\frac{uv}{-1+u+v+\big[(1-u)^m+(1-v)^m\big]^{1/m}}
\end{equation} via the formula (\ref{extremeC3}). Consequently, $C_m(F_1,F_2)$ is a bi-free extreme value distribution function whenever $F_1$ and $F_2$ are of free extreme type. The copula $C_m$ is in the domain of attraction of the extreme-value copula
\[\exp\left\{-\big[(-\log u)^m+(-\log v)^m\big]^{1/m}\right\},\] known as the \emph{Gumbel-Hougaard copula}.}
\end{exam}

\begin{exam} \emph{With the content in the previous example, $A_m$ tends to $A_\infty(t)=1\vee(1-t)$ as $m\to\infty$. Using (\ref{extremeC3}) with the dependence function $A_\infty$, the comonotone copula is obtained, which is also the extreme-value copula determined by $A_\infty$. In other words, the comonotone copula can be used to generate not only extreme value distribution functions but also bi-free extreme value ones.}
\end{exam}

\begin{exam} \emph{(Non-differentiable Model). Let $\theta,\phi\in[0,1]$. Elementary computations show that the copula
\[\frac{uv}{1-\min\big\{\theta(1-u),\phi(1-v)\big\}}\] lies in the domain of attraction of the \emph{Marshall-Olkin copula}
\[uv\min\big\{u^{-\theta},v^{-\phi}\big\},\] which has the dependence function $1-\min\{\theta t,\phi(1-t)\}$. Once again, these non-differentiable copulas can be used to respectively generate bi-free extreme and extreme value distribution functions.}
\end{exam}

\section*{Acknowledgments} The first author
was supported by a grant from the Ministry of Science and Technology in Taiwan MOST 106-2628-M-110-002-MY4. The second author was supported by the
NSERC Canada Discovery Grant RGPIN-2016-03796.

\end{document}